\documentclass[onefignum,onetabnum,onealgnum]{siamart220329}
\pdfoutput=1

\usepackage{lipsum}
\usepackage{amsfonts}
\usepackage{graphicx}
\usepackage{amsopn}
\DeclareMathOperator{\diag}{diag}
\usepackage{mathtools}
\usepackage{paralist}

\usepackage{color}
\usepackage[acronym,shortcuts]{glossaries}

\usepackage{pgf}
\usepackage{algpseudocode} % for 'algorithmicx' which is an improved version of algorithmic
\usepackage{csquotes}

\usepackage{tikz}
\usetikzlibrary[arrows, decorations.pathmorphing, backgrounds, patterns, positioning, fit, petri]

\usepackage{siunitx}
\usepackage[symbol]{footmisc}

\usepackage{booktabs}

\usepackage{makecell}

\newtheorem{remark}{Remark}

\algrenewcommand{\algorithmiccomment}[2][.6\linewidth]{\leavevmode\hfill\makebox[#1][l]{\textcolor{darkgray}{$\triangleright$ #2}}}
\algnewcommand{\algorithmicgoto}{\textbf{goto}}%!TEX encoding = UTF-8 Unicode
\algnewcommand{\Goto}[1]{\algorithmicgoto~\ref{#1}}%
\algnewcommand\algorithmicvariables{\textbf{Local variables:}}%
\algnewcommand\Variables{\item[\algorithmicvariables]}%

\algblock{ParFor}{EndParFor}
\algnewcommand\algorithmicparfor{\textbf{parallel for}}
\algnewcommand\algorithmicpardo{\textbf{do}}
\algnewcommand\algorithmicendparfor{\textbf{end\ parallel for}}
\algrenewtext{ParFor}[1]{\algorithmicparfor\ #1\ \algorithmicpardo}
\algrenewtext{EndParFor}{\algorithmicendparfor}

\algnewcommand{\LineComment}[1]{\vspace{3pt}\Statex \(\triangleright\)\(\triangleright\) \textcolor{blue}{#1} \vspace{3pt}}
\newcommand{\COMMENT}[2][.47\linewidth]{%
    \leavevmode\hfill\makebox[#1][l]{$\triangleright$~#2}}

\newcommand{\ir}{$\mathcal{IR}$}
\newcommand{\V}{$\mathcal{V}$}
\newcommand{\fmg}{$\mathcal{FMG}$}
\newcommand{\s}{$\mathcal{S}$}

\newcommand{\exact}[1]{{#1}^\ast}

\newcommand{\bfp}{\mathbb{B}}
\newcommand{\bfpm}[1]{#1_m} % block-float "mantissa vector"
\newcommand{\bfpe}[1]{#1_e} % block-float exponent
\newcommand{\bfpq}[1]{#1_w} % block-float mantissa length
\newcommand{\bfpd}[1]{#1_d} % block-float dimension(s)
\newcommand{\bfpt}[1]{(\bfpe{#1}, \bfpm{#1}, \bfpq{#1}, \bfpd{#1})} % block-float tuple
\newcommand{\bfpequal}{\sim}

\newcommand{\LeftShift}{\ll}
\newcommand{\RightShift}{\gg}

\newcommand{\real}{\mathbb{R}}
\newcommand{\unsigned}{\mathbb{Z}_{\geq 0}}
\newcommand{\unsignedpositive}{\mathbb{Z}_{> 0}}
\newcommand{\signed}{\mathbb{Z}}

\newcommand{\twoscomp}[1]{\mathbb{X}_{#1}}     % integer range for bit-width w (twos complement)
\newcommand{\tmp}[1]{{#1}_\textrm{tmp}}                 % temporary variable
\newcommand{\out}[1]{{#1}_\textrm{out}}                 % output variable

\newcommand{\ewe}{\boldsymbol{\varepsilon}}
\newcommand{\ewed}{\boldsymbol{\dot{\varepsilon}}}
\newcommand{\eweb}{\boldsymbol{\bar{\varepsilon}}}
\newcommand{\eweq}{\boldsymbol{\check{\varepsilon}}}

\newcommand{\wewe}{w}
\newcommand{\wewed}{\dot{w}}

\newcommand{\weweq}{\check{w}}

\definecolor{ColorNeutralBits}{RGB}{35, 87, 137}
\definecolor{ColorEnoughBits}{RGB}{59, 178, 115}
\definecolor{ColorOverflow}{RGB}{239, 48, 84}
\definecolor{ColorUnderflow}{RGB}{255, 111, 89}

\definecolor{lo}{HTML}{33A02C}
\definecolor{med}{HTML}{1F78B4}
\definecolor{hi}{HTML}{E31A1C}

\definecolor{lo_light}{HTML}{B2DF8A}
\definecolor{med_light}{HTML}{A6CEE3}
\definecolor{hi_light}{HTML}{FB9A99}
\definecolor{quant}{HTML}{FF7F00}

\newacronym{bfp}{BFP}{block floating point}
\newacronym{blas}{BLAS}{Basic Linear Algebra Subprograms}
\newacronym{fmg}{FMG}{full multigrid}
\newacronym{msb}{MSB}{most significant bit}
\newacronym{lsb}{LSB}{least significant bit}
\newacronym{pde}{PDE}{partial differential equation}
\newacronym{spd}{SPD}{symmetric positive definite}

\headers{Multigrid methods using Block Floating Point Arithmetic}{N. Kohl, S. F. McCormick, and R. Tamstorf}

\title{Multigrid methods \\ using Block Floating Point Arithmetic%
}

\author{Nils Kohl\thanks{Friedrich-Alexander-Universität Erlangen-N{\"u}rnberg, Germany (\email{nils.kohl@fau.de}).}
\and Stephen F. McCormick\thanks{University of Colorado at Boulder, Boulder, CO (\email{stephen.mccormick@colorado.edu}).}
\and Rasmus Tamstorf\thanks{Walt Disney Animation Studios, Burbank, CA (\email{rt@acm.org}).}}

\begin{document}

    \maketitle

% ---------------------------------------------
% ---------------------------------------------

    \begin{abstract}
      \Gls*{bfp} arithmetic is currently seeing a resurgence in interest because it requires less power, less chip area, and
      is less complicated to implement in hardware than
      standard floating point arithmetic. This paper explores the application of \gls*{bfp} to
      mixed- and progressive-precision multigrid methods, enabling
      the solution of linear elliptic \glspl*{pde} in energy- and hardware-efficient \emph{integer arithmetic}.
      While most existing applications of \gls*{bfp} arithmetic tend to use small block sizes, the block size here is chosen to be maximal such that
        matrices and vectors share a \emph{single} exponent for all entries. This is sometimes also referred to as a scaled fixed-point format.
        We provide algorithms for \acs*{blas}-like routines for \gls*{bfp} arithmetic that ensure exact vector-vector
        and matrix-vector computations up to a specified precision.
        Using these algorithms, we study the asymptotic precision requirements to achieve
        discretization-error-accuracy.
        We demonstrate that some computations can be performed using as little as 4-bit integers, while the 
        number of bits required to attain a certain target accuracy is similar to that of standard floating point
        arithmetic.
        Finally, we present a heuristic for full multigrid in \gls*{bfp} arithmetic based on saturation and truncation
        that still achieves
        discretization-error-accuracy without the need for expensive normalization steps of intermediate results.
    \end{abstract}

    \begin{keywords}
        % 7. Keywords that describe the paper
        Block floating point, fixed point, mixed precision, multigrid
    \end{keywords}


    \section{Introduction}\label{sec:intro}

    Floating point arithmetic is used to perform almost all scientific computations.
    At the same time, it is well known that integer arithmetic, where applicable, is less complicated to implement in hardware than
    standard floating point arithmetic~\cite{Behrooz:2010} and it is generally more energy efficient. As an example, the actual arithmetic associated with 32 bit integer based addition requires roughly an order of magnitude less energy than the corresponding floating point operation, \cite{Horowitz2014,Jouppi2021}.
    For this reason, \emph{fixed point} formats are typically preferred in embedded computing where resources are limited.
    Unfortunately, the inherent range limitation of fixed point formats renders them difficult to use for the numerical approximation of partial differential equations 
    (\glspl*{pde}).
    A compromise is to use a \glsfirst*{bfp} format:
    a block of fixed point mantissas along with a \emph{shared exponent}~\cite{Wilkinson1963}.
    In this way, the range of representable numbers in \gls*{bfp} formats can be adapted dynamically, while all computations are
    still performed in pure integer arithmetic. While the notion of 
    \gls*{bfp} and fixed point formats goes back quite far in the history of computing, it has recently gained renewed popularity for neural network
    training, e.g.,~\cite{Koster2017,Drumond2018,Lian2019,Rouhani2020,Dai2021,QianZhang:2022:FASTDNNTraining,Noh2022a,Noh2022b,Basumallik2022}. The actual cost of arithmetic is typically dwarfed by the cost of memory access, but \cite{Noh2022a} shows that it is possible to increase the overall energy efficiency by an order of magnitude when using \gls*{bfp} compared to using the mixed FP16/FP32 arithmetic in Nvidia's tensor cores. 
   
    In this paper, we study the solution of linear systems arising from the discretization of elliptic
    \glspl*{pde} in \gls*{bfp} arithmetic using mixed- and progressive-precision multigrid methods.
    As in~\cite{McCormick:2021:AlgebraicErrorAnalysis,Tamstorf:2021:DiscretizationErrorAccurateMixedPrecisionMultigrid},
    we are interested in the asymptotically optimal choice of precisions that guarantees
    discretization-error-accurate solutions.
    We design algorithms for matrix-vector and vector-vector operations in \gls*{bfp}-arithmetic that ensure efficient
    and exact computations up to a specified target precision, and emphasize that all computations are performed in
    two's complement integer arithmetic.

    To ensure exact computations, we leverage the fact that \gls*{bfp} enables the computation of the exact inner product between two vectors at a reasonable cost \cite{Drumond2018}.
The exact dot product can also be computed for standard floating point numbers using the method proposed by Kulisch \cite{Kulisch2011b}, and this method is used by posits to implement the socalled ``quire'', \cite{Gustafson2017}. However, in general it requires a very large accumulator with more than 4,000 bits for double precision floating point numbers and more than 65,000 bits for quad-precision numbers. For block floating point numbers, the bulk of the computation can be done in fixed point arithmetic where all the elements are stored in the same format. The size of the accumulator for the result of the exact \gls*{bfp} dot product therefore only grows logarithmically with the number of vector entries. Thus, the technique in \cite{Boldo2020} can be used to compute an exact dot product for high-precision inputs using on the order of $100-200$ bits or less for most practical cases.

    The outline of the remaining parts of the paper is as follows:
    \Cref{sec:bfp} introduces the \gls*{bfp}-format and its most relevant properties.
    \Cref{sec:bfp-quant-error} analyses the relative energy error induced by \gls*{bfp}-quantization.
    \Cref{sec:bfp-blas} develops
    \gls*{bfp}-specific algorithms for mixed-precision matrix-vector and vector-vector operations.
    \Cref{sec:mg} summarizes the results
    of~\cite{McCormick:2021:AlgebraicErrorAnalysis,Tamstorf:2021:DiscretizationErrorAccurateMixedPrecisionMultigrid}
    and extends the mixed-precision multigrid method defined therein to \gls*{bfp}-arithmetic.
    \Cref{sec:results} provides a numerical study of the precision requirements for the \gls*{bfp}-multigrid solver
    in order to achieve discretization-error-accurate approximations in the energy norm for two model problems. We end in \Cref{sec:conclusion} with concluding remarks.

%%%%%%%%%%%%%%%%%%%%%%%%%%%%%%%%%%%%%%%%%%%%%%%%%%%%%%%%%%%%%%%%%%%%%%%%%%%%%%%%%%%%%%%%%%%%%%%%%%%%%%%%%%%%%%%%%%%%%%%%
%%%%%%%%%%%%%%%%%%%%%%%%%%%%%%%%%%%%%%%%%%%%%%%%%%%%%%%%%%%%%%%%%%%%%%%%%%%%%%%%%%%%%%%%%%%%%%%%%%%%%%%%%%%%%%%%%%%%%%%%
%%%%%%%%%%%%%%%%%%%%%%%%%%%%%%%%%%%%%%%%%%%%%%%%%%%%%%%%%%%%%%%%%%%%%%%%%%%%%%%%%%%%%%%%%%%%%%%%%%%%%%%%%%%%%%%%%%%%%%%%

    \section{Block floating point arithmetic}\label{sec:bfp}

    We define \gls*{bfp} numbers by a block of integers (also referred to as mantissas)
    equipped with a shared factor that is an integer power of two.
    All integers of a block have the same bit-width, and we assume standard two's complement
    representation.
    We denote the set of two's complement integers with bit-width $w \in \unsignedpositive$ as
    $\twoscomp{w} \coloneqq [-2^{w-1}, 2^{w-1} - 1] \cap \signed$.
    The block will typically be a vector or matrix, but the concept generalizes to any type of tensor or irregular
    structures such as sparse matrix formats.
    We denote the set of \gls*{bfp} numbers by $\bfp$ and write each element as a tuple
    \begin{equation}
        \bfpt{x} \in \bfp.
    \end{equation}
    In this notation, $\bfpe{x}\in\signed$ denotes the shared exponent, and $\bfpq{x}\in\unsignedpositive$ denotes the fixed
    bit-width of the mantissas.
    With a slight abuse of notation, we let $\bfpm{x} \in \twoscomp{\bfpq{x}}^{\bfpd{x}}$ denote the block of
    mantissas, where $\bfpd{x}$ represents the layout of the elements in the block.
    In the case of an $n$-dimensional vector, we let $\bfpd{x}=n$, and in the case of an $n\times n$ matrix, we
    write $\bfpd{x} = n\times n$.
    We allow the exponent $\bfpe{x}$ to be chosen arbitrarily because it is shared over the entire block and its
    storage cost is negligible in practice. (A 64-bit integer exponent is likely more than sufficient for most practical use
    cases.)
    We use the notation $x \bfpequal \bfpt{x}$ to denote that $x$ is the block of rational numbers
    $x = 2^{\bfpe{x}} \cdot \bfpm{x}$, and the shorthand $\bfp^{d} = \{ \bfpt{x} \in \bfp : \bfpd{x} = d \}$.

%    \subsection{Normalization and precision}

    Generally, multiple equivalent representations of $x$ correspond to different choices of $\bfpe{x}$.
    We call the representation $\bfpt{x}$ of $x \neq \mathbf{0}$ \emph{normalized} if $\bfpe{x}$ is
    minimal (possibly negative).
    In the following, $x\in\bfp$ is used to refer to both the tuple containing the representation for the
    \gls*{bfp} numbers and the represented numbers.

    Relevant properties of a \gls*{bfp} format can be derived from the quantities $\bfpe{x}$ and $\bfpq{x}$.
    As an example, the range of all entries of $x$ is
    $\left[-2^{\bfpq{x}-1} \cdot 2^{\bfpe{x}}, (2^{\bfpq{x}-1} - 1) \cdot 2^{\bfpe{x}} \right]$, and all the representable
    numbers are equidistantly separated with distance $2^{\bfpe{x}}$.
    The precision of a normalized \gls*{bfp} tensor with a fixed mantissa width $\bfpq{x}$
    therefore depends on the value of the entry that has the largest magnitude.
    To relate floating point precision (i.e.,~unit roundoff) to the \gls*{bfp} context, we refer in this paper to the
    precision of a \gls*{bfp} format as $\varepsilon = 2^{-(\bfpq{x} - 1)}$.
    This is the spacing between two adjacent values in a normalized \gls*{bfp} tensor with entries in the range $[-1, 1 - \varepsilon]$.

    We use standard two's complement integer arithmetic for addition, subtraction, and multiplication, as well as arithmetic left- ($\LeftShift$) and
    right-shifts ($\RightShift$). Additionally, we define the operation $\text{decr}(b, \cdot)$ that truncates the $b \geq 0$
    leftmost bits of a two's complement integer (corresponding to casting to a narrower integer type),
    and the operation $\text{incr}(b, \cdot)$ that prepends $b$ bits that all have the value of the \gls*{msb} to the
    left (corresponding to casting to a wider type).
    Truncation of the rightmost bits is realized via arithmetic right-shifts, which implies rounding
    towards negative infinity for signed two's complement integers.

    \begin{remark}[Block size]
        Throughout this paper, we consider the extreme case that each vector or matrix is represented by a single block.
        At the other extreme, setting the block size $\bfpd{x}$ to 1 is equivalent to using standard floating point arithmetic.
        From an implementation standpoint, these extreme cases (block sizes 1 and maximal) are special,
        as block-boundaries can be ignored.
        However, there may be practical reasons to split up a vector or matrix into multiple blocks.
        As an example, hardware components may be specialized to perform optimized arithmetic on relatively small
        block sizes (e.g., on the order $10-100$ entries)~\cite{Rouhani2020,Dai2021}.
        That being said, choosing the block size maximally has to be the worst case in terms of quantization.
        Thus, the results presented herein are expected to extend easily to smaller block sizes.
    \end{remark}

    \begin{remark}[\gls*{bfp} dot products]
        The energy saving potential of \gls*{bfp} arithmetic compared to standard floating point arithmetic is rooted in the simplification of the dot product.
        This carries over to matrix-vector and matrix-matrix multiplications, which conceptually are just consecutive dot products.
        When solving sparse linear systems using iterative solvers, the critical
        steps to performance are sparse matrix-vector multiplications, which makes
        an efficient dot product particularly beneficial.

        To understand the origin of the complexity reduction, consider the addition of two floating point numbers.
        Before the mantissas can be added, they have to be aligned by right-shifting one of them in order to ensure that the exponents
        are equal. After the addition, the result has to be normalized, which requires another shift operation.
        This is somewhat simplified, but illustrates the complexity of a seemingly simple operation.
        A dot product of two vectors $a, b \in \real^n$ with at most $m_A$ non-zero elements per vector
        %($a$ representing the non-zeros of a row of $A$ and $b$ the corresponding entries of $x$)
        requires $m_A - 1$ additions and, therefore, if
        performed in floating point arithmetic, $2 (m_A - 1)$ arithmetic shifts. 
        This is different in \gls*{bfp} arithmetic: assuming that the block sizes are maximal, the alignment
        step is not necessary at all, since the terms in the sum in the dot product all share the
        same exponent. (This exponent is computed by adding the block-exponents of $a$ and $b$.)
        Furthermore, normalization is only necessary after summing up the result in a sufficiently large register.
        In other words, only a single arithmetic shift is necessary in \gls*{bfp} arithmetic compared to
        $2 (m_A - 1)$ arithmetic shifts in floating point arithmetic. 

        A key assumption here is that the accumulator is sufficiently large to hold the sum without overflow. Since
        the number of non-zero terms in the sum is assumed to be at most $m_A$, it follows that the size of the accumulator
        grows with $\log_2(m_A)$. In practice, $m_A$ is typically at most in the hundreds, so only a relatively few additional
        bits are required to compute the exact dot product.
        With suitable rounding, this number could possibly be reduced.
    \end{remark}

    \section{Relative BFP-quantization error}\label{sec:bfp-quant-error}
    
    A critical issue concerning the practical use of \gls*{bfp} arithmetic is its effect on accuracy. The empirical observations in~\cite{McCormick:2021:AlgebraicErrorAnalysis,Tamstorf:2021:DiscretizationErrorAccurateMixedPrecisionMultigrid} suggest that the relative fixed-point quantization error  is $\mathcal{O}(\kappa^{\frac{1}{2}} \ewe)$, and our experience indicates that \gls*{bfp} exhibits the same order. The aim in this section is to develop theoretical results that shed more light on this issue. In particular, we provide an abstract bound that suggests that \gls*{bfp} quantization might incur a slightly higher order of error. We then argue that this bound might be pessimistic in that it does not fully take finite-precision into account.

To be specific, let $A$ denote the \gls*{spd} system matrix of an elliptic \gls*{pde} discretized by standard finite elements on a uniform $n \times n$ grid in the unit square. Assume for simplicity that $A$ is scaled so that its minimal eigenvalue is $\mathcal{O}(1)$. With $u$ an infinite-precision vector in $\real^{n^2}$ having unit infinity norm $\|u\|_\infty$, write $u = v + \ewe z$, where $v$ results from \gls*{bfp} quantization in $\ewe$ precision and $\ewe z$ is the quantization error. Note that $\|u\|_\infty = \|v\|_\infty = 1$ and $\|z\|_\infty \le 1$, where $\|\cdot\|_\infty$ denotes the infinity norm. The goal is then to bound the relative \gls*{bfp} quantization error $\mathcal{E}(v, z) \coloneqq \frac{\ewe \|z\|_A}{\|v\|_A}$, where $\|\cdot\|_A = \|A^\frac{1}{2}\cdot\|$ denotes the energy norm written in terms of the Euclidean norm $\|\cdot\|$.

We can usually choose $z$ so that it has maximum order, that is,
\begin{equation}
\|z\|_A = \mathcal{O}(\kappa^{\frac{1}{2}}n),
\label{osc}
\end{equation}
while preserving the property that the \gls*{bfp} quantization of $u$ is $v$. For example, with the five-point discrete 2D Poisson equation, we could choose $z$ to alternate between $1$ and $0$ in a checkerboard fashion. This choice means that \gls*{bfp} truncation of $u$ does indeed result in $v$ (as might not be the case with negative values of $v$) and it also assures that (\ref{osc}) holds: $z$ is oscillatory in that $\|z\|_A = \mathcal{O}(\kappa^{\frac{1}{2}}) \|z\|$ and it possesses enough $1$'s to make $\|z\| = \mathcal{O}(n)$. To bound $\mathcal{E}(v, z)$ in this case, we therefore need only find a lower bound for $\|v\|_A$.  

We have not been able to establish a sharp lower bound for $\|v\|_A$ theoretically because it ostensibly requires the discrete optimization of $\|v\|_A$ over the space of vectors of unit infinity norm that are represented exactly in $\ewe$ precision. We can, however, obtain a potentially loose lower bound by ignoring the finite-precision restriction and exploiting the fact that the minimum value of $\|v\|_A$ is the inverse of the square root of the maximum of the diagonal entries of $A^{-1}$. (See Appendix~\ref{disc-harm}.) Our numerical estimates of $A^{-1}$ for the model 2D Poisson problem for $n \in \{1, 2, \dots, 100\}$ indicate that the minimum value of $\|v\|_A$ is bounded below by a constant times $n^{0.994}$, suggesting that $\mathcal{E}(v, z)$ might grow slightly faster than $\mathcal{O}(\kappa^{\frac{1}{2}} \ewe)$ (by a factor of $n^{0.006}$).

\begin{figure}
        \centering
        \resizebox*{\textwidth}{!}{%
            \input{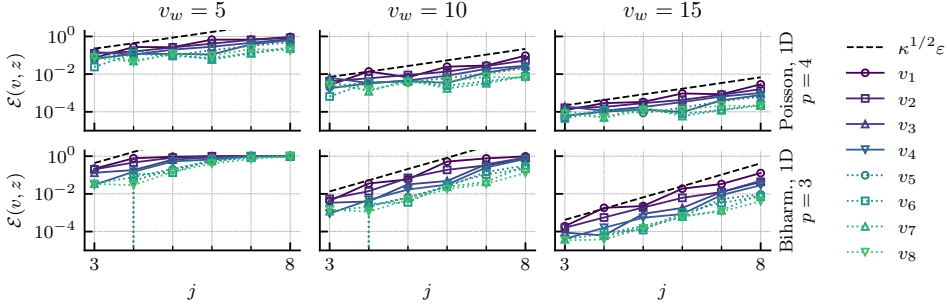}%
        }
        \vspace{-24pt}
        \caption{Relative energy errors $\mathcal{E}(v, z) \coloneqq \ewe \|z\|_A / \|v\|_A$ after quantization of
        $8$ eigenvectors $u_i$ of $A$ for the model problem discussed in \Cref{sec:results}. The precision is indicated by the mantissa width
        $\bfpq{v} \in \{5, 10, 15\}$.
        The x-axis shows the refinement level, $j$, with corresponding mesh sizes $h = 2^{-j}$.}
        \label{fig:bfp-relative-quant-errors}
    \end{figure}

The slightly larger bound requires $u$ to be very smooth while the part $\ewe z$ that is truncated away is oscillatory. While we do not know how likely this is, we have not experienced an error growth that is larger than $\mathcal{O}(\kappa^{\frac{1}{2}} \ewe)$. Just for illustration, \cref{fig:bfp-relative-quant-errors} shows $\kappa^{\frac{1}{2}} \ewe$ and
    the relative energy error after quantization to \gls*{bfp} over refinement for different mantissa widths $\bfpq{v}$.
    The precision $\ewe$ is computed as described in \Cref{sec:bfp}.
    We chose $v$ as the quantized eight eigenvectors $v_i$ of $A$ that belong to the eight smallest eigenvalues
    $\lambda_i,\ i = 1, \dots, 8$, since those yielded the largest relative errors. (Note that the errors are getting
    smaller as $i$ increases.) The model problems are discussed in \Cref{sec:results}. 
    
    In any case, a possible reason for the discrepancy between theory and our experience is that the theory is only an upper bound. Indeed, since $v$ has not been restricted to $\ewe$-precision in this theoretical bound, it may be an overestimate due to the theoretical minimum being taken over a wider set. 

    \section{BLAS operations using BFP}\label{sec:bfp-blas}

    The majority of iterative solvers for linear systems of equations is composed of vector-vector and matrix-vector
    operations.
    In our case, most computations can be written in the form
    \begin{equation}\label{eq:gemv}
        z \gets \alpha A x + \beta y,
    \end{equation}
    with $y, z \in \mathbb{R}^n$, $x \in \mathbb{R}^m$, $A \in \mathbb{R}^{n \times m}$, $\alpha, \beta \in \mathbb{R}$.
    This operation is also known as the generalized matrix-vector multiplication (\texttt{gemv}) from the \gls*{blas}
    specification~\cite{BasicLinearAlgebraSubprogramsTechnicalForum:2002:BasicLinearAlgebra}.
    It can be simplified in some cases, such as when $A = I$ (\texttt{axpby}), $y = \mathbf{0}$, or
    $\alpha = 1$.

    One contribution of this paper is the realization of a framework for such vector operations that is suited for mixed-precision computations
     using \gls*{bfp} arithmetic with a
    focus on the requirements for multigrid methods.
    We assume that
    all inputs are given in a normalized \gls*{bfp} format with arbitrary, and possibly different, mantissa
    widths and block-exponents.
    The desired width of the mantissa of the result, denoted $\out{w}$, is assumed to be specified as an input parameter, and 
    the result is required to be normalized and computed \emph{exactly} up to chosen size of the mantissa.

    The normalization of the result is not straightforward, since we generally do not know the largest \gls*{msb}
    index of the block-mantissa before the entire vector has been computed.
    Two naive implementations come to mind:
    \begin{itemize}
        \item[(a)] (\emph{computationally efficient})
        The entire result is computed exactly and \emph{stored} in a temporary unnormalized \gls*{bfp} vector using a sufficiently wide mantissa.
        The \gls*{msb} index is tracked.
        In a second step, the block-mantissa is shifted to normalize the result, and then quantized to $\out{w}$.

        \item[(b)] (\emph{memory efficient})
        The result is first computed exactly element by element only to track the \gls*{msb} index for the entire block.
        The actual result is discarded after the computation of each element.
        In a second step, the result is \emph{recomputed} element by element with each element being shifted to ensure that the overall block is normalized. The result is also quantized to $\out{w}$ bits.
    \end{itemize}

    Both approaches have disadvantages.
    In (a), the entire temporary BFP vector carrying the exact result must be allocated.
    Depending on the precision and block-exponents of the input variables, the amount of memory required can be huge.
    In (b), the exact result is not stored at all.
    Only the quantized result is stored, but this roughly doubles the computational effort.
    We suggest a compromise of these two approaches, where the temporary result is neither stored exactly nor
    discarded but rather quantized to some intermediate width $\tmp{w} \geq \out{w}$.
    The idea is to choose the width large enough so that at least $\out{w}$ ``good'' bits are left for the largest element after normalization of the entire block,
    but still small enough so that much less memory is allocated than would be required for the exact result.

    This algorithm requires a (preferably sharp) upper bound $\gamma \geq \|z\|_\infty$ and a
    number of bits $\tmp{w}$ to prescribe a ``window'' in the mantissa of the exact result that ideally contains the
    most significant $\out{w}$ bits.
    When the result is computed, only those $\tmp{w}$ bits that make up that window are stored.
    Everything beyond the window is discarded.

    In the following, indices are used to refer to individual bits of a bit string.
    The rightmost bit is assigned index $0$, with incrementing indices to the left.
    To refer to a substring of a bit string, we use
    $\mu$ (\underline{m}u) to refer to the index left to the \glsfirst*{msb} (leftmost bit) of that substring,
    and $\lambda$ (\underline{l}ambda) to refer to the index of its \gls*{lsb} (rightmost bit).
    For convenience, $\mu$ is set to the bit-index to the left of the actual \gls*{msb} of the substring,
    so that the width of the substring is $\mu - \lambda$.
    For example, given a bit string $\texttt{00010010}$, its substring $\texttt{1001}$ is indicated by $\mu = 5$, and
    $\lambda = 1$.
    Its width is $\mu - \lambda = 4$.

    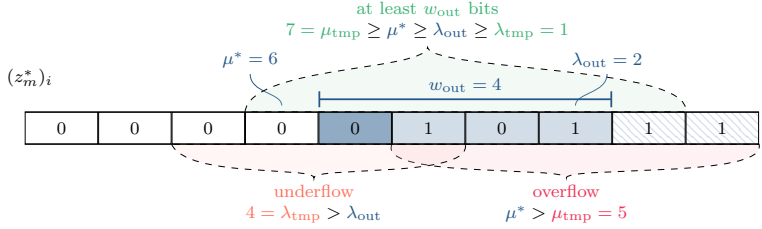
\begin{figure}[t]
        \centering
        \resizebox*{.8\textwidth}{!}{%
            \resizebox{0.9\textwidth}{!}{%
    \begin{tabular}{@{}l@{}}
        \begin{tikzpicture}[node distance=0pt,
            box/.style={draw, minimum size=0.5cm, inner xsep=0.5cm, inner ysep=0.3},
            emptybox/.style={draw, dashed, minimum size=1cm, inner sep=0.5cm},
            value/.style={yshift=-1.5cm}]

            \footnotesize

            \node[box]                                                              (b9)               {0};
            \node[box]                                                              (b8) [right=of b9] {0};
            \node[box]                                                              (b7) [right=of b8] {0};
            \node[box]                                                              (b6) [right=of b7] {0};
            \node[box, fill=ColorNeutralBits!50]                                    (b5) [right=of b6] {0};
            \node[box, fill=ColorNeutralBits!20]                                    (b4) [right=of b5] {1};
            \node[box, fill=ColorNeutralBits!20]                                    (b3) [right=of b4] {0};
            \node[box, fill=ColorNeutralBits!20]                                    (b2) [right=of b3] {1};
            \node[box, pattern=north west lines, pattern color=ColorNeutralBits!20] (b1) [right=of b2] {1};
            \node[box, pattern=north west lines, pattern color=ColorNeutralBits!20] (b0) [right=of b1] {1};

%                    \node[yshift=2.8cm] [above of=b9] {$9$};
%                    \node[yshift=2.8cm] [above of=b8] {$8$};
%                    \node[yshift=2.8cm] [above of=b7] {$7$};
%                    \node[yshift=2.8cm] [above of=b6] {$6$};
%                    \node[yshift=2.8cm] [above of=b5] {$5$};
%                    \node[yshift=2.8cm] [above of=b4] {$4$};
%                    \node[yshift=2.8cm] [above of=b3] {$3$};
%                    \node[yshift=2.8cm] [above of=b2] {$2$};
%                    \node[yshift=2.8cm] [above of=b1] {$1$};
%                    \node[yshift=2.8cm] [above of=b0] {$0$};

            \node[yshift=0.8cm, xshift=-0.5cm] [above of=b9] {$(\bfpm{\exact{z}})_i$};
%
%                    \draw [|-|, color=black, thick] ([yshift=1.4cm]b9.north west) to
%                    node[color=black, midway, above] {$\bfpq{\exact{r}} = 10$} ([yshift=1.4cm]b0.north east);
%
%                    \draw [|-|, color=black, thick] ([yshift=0.8cm]b5.north west) to
%                    node[color=black, midway, above] {$\out{w} = 4$} ([yshift=0.8cm]b2.north east);
%
%                    \draw [|-|, color=black, thick] ([yshift=0.2cm]b6.north west) to
%                    node[color=black, midway, above] {$\temp{w} = 6$} ([yshift=0.2cm]b1.north east);

            \node[text=ColorNeutralBits, align=center, xshift=-0.5cm, yshift=0.6cm] (nb6) [above=of b6] {$\exact{\mu} = 6$};
            \draw [shorten >=0.1cm, color=ColorNeutralBits, in=90, out=270]     (nb6.south) to (b6.north);

            \node[text=ColorNeutralBits, align=center, xshift=0.5cm, yshift=0.6cm] (nb2) [above=of b2] {$\out{\lambda} =
            2$};
            \draw [shorten >=0.1cm, color=ColorNeutralBits, in=90, out=270]     (nb2.south) to (b2.north);

%                    \node[value, text=green, align=center] (nb7) [below of=b7] {$\tmp{\mu} = 7$};
%                    \draw [->, shorten >=0.1cm, shorten <=0.1cm, color=green]     (nb7.north) to (b7.south);
%
%                    \node[value, text=green, align=center] (nb1) [below of=b1] {$\tmp{\lambda} = 1$};
%                    \draw [->, shorten >=0.1cm, shorten <=0.1cm, color=green]     (nb1.north) to (b1.south);

            \node[align=center] at (5.8, 1.7) (success) {\textcolor{ColorEnoughBits}{at least $\out{w}$ bits} \\
            $\textcolor{ColorEnoughBits}{7 = \tmp{\mu}} \geq \textcolor{ColorNeutralBits}{\exact{\mu}} \geq
            \textcolor{ColorNeutralBits}{\out{\lambda}} \geq \textcolor{ColorEnoughBits}{\tmp{\lambda} = 1}$};

            \draw[fill=ColorEnoughBits!30, fill opacity=0.2, dashed] (b7.north east) .. controls (3.0, 1.0) and (5.8, 0.9) .. (success.south)
            -- (success.south) .. controls (5.8, 0.9) and (10, 1.0) .. (b1.north east)
            -- cycle;

            \node[align=center] at (8, -1.2) (failoverflow) {\textcolor{ColorOverflow}{overflow} \\
            $\textcolor{ColorNeutralBits}{\exact{\mu}} > \textcolor{ColorOverflow}{\tmp{\mu} = 5}$};

            \draw[fill=ColorOverflow!30, fill opacity=0.2, dashed] (b5.south east) .. controls (5.0, -0.8) and (8, -0
            .6) .. (failoverflow.north)
            -- (failoverflow.north) .. controls (8, -0.6) and (11.2, -0.8) .. (b0.south east)
            -- cycle;

            \node[align=center] at (4, -1.2) (failbits) {\textcolor{ColorUnderflow}{underflow} \\
            $\textcolor{ColorUnderflow}{4 = \tmp{\lambda}} > \textcolor{ColorNeutralBits}{\out{\lambda}}$};

            \draw[fill=ColorUnderflow!30, fill opacity=0.2, dashed] (b8.south east) .. controls (2, -0.8) and (4, -0.6)
            .. (failbits.north)
            -- (failbits.north) .. controls (4, -0.6) and (6.5, -0.8) .. (b4.south east)
            -- cycle;

            \draw [] (b9.north west) to (b0.north east);
            \draw [] (b9.south west) to (b0.south east);

            \draw [|-|, color=ColorNeutralBits, thick] ([yshift=0.2cm]b5.north west) to
            node[color=black, midway, above] {$\textcolor{ColorNeutralBits}{\out{w} = 4}$} ([yshift=0.2cm]b2.north east);

        \end{tikzpicture}
    \end{tabular}
}%
        }
        \caption{Illustration of the bit indices in \cref{alg:qcomp}.
        The bit string $(\bfpm{\exact{z}})_i$ shown here represents one element of the block-mantissa $\bfpm{\exact{z}}$ of
        the exact result $\exact{z} \in \bfp^n$. We assume for simplicity that this element has the maximum \gls*{msb} index
        of all elements in $\bfpm{\exact{z}}$.
        The actual \gls*{msb} and \gls*{lsb} indices of the relevant portion (filled boxes) of $(\bfpm{\exact{z}})_i$
            are $\textcolor{ColorNeutralBits}{\exact{\mu}}$ and $\textcolor{ColorNeutralBits}{\out{\lambda}}$.
            We show three example bit-``windows'' that are \textcolor{ColorEnoughBits}{sufficient} (top),
            or result in \textcolor{ColorUnderflow}{underflow} (bottom left), or \textcolor{ColorOverflow}{overflow}
            (bottom right).
            Only if the estimated \gls*{msb} and \gls*{lsb} indices $\textcolor{ColorEnoughBits}{\tmp{\mu}}$ and
            $\textcolor{ColorEnoughBits}{\tmp{\lambda}}$ (example choices are displayed in the figure) that are
            derived from $\gamma$ and $\tmp{w}$ fulfill
            $\textcolor{ColorEnoughBits}{\tmp{\mu}} \geq \textcolor{ColorNeutralBits}{\exact{\mu}} \geq
            \textcolor{ColorNeutralBits}{\out{\lambda}} \geq \textcolor{ColorEnoughBits}{\tmp{\lambda}}$,
            are all
            relevant bits (filled boxes) captured by the window. Otherwise, we either obtain underflow or overflow
            (illustrated by example choices for $\textcolor{ColorUnderflow}{\tmp{\lambda}}$ and
            $\textcolor{ColorOverflow}{\tmp{\mu}}$).}
        \label{fig:windows-merged}
    \end{figure}

    %The hope is that we can roughly estimate a rather sharp upper bound to the result.
    %For example this should be possible for the residual computation in \ir~ using the previous residual and an
    %estimate of the convergence rate of the inner solver.

    After looping through the entire result vector, we know the \gls*{msb} index $\exact{\mu}$ of the exact mantissa
    (we denote variables corresponding to exact quantities using asterisks).
    Therefore, we can determine whether at least $\out{w}$ ``good'' bits have been collected.
    If that is the case, the temporary result is shifted and quantized to $\out{w}$ bits.
    We slightly abuse the terms overflow and underflow to refer to the cases where we either do not capture the
    \gls*{msb} or capture less than $\out{w}$ relevant bits of the mantissa entry with the maximum \gls*{msb} index.
    In such cases, the result has to be recomputed.
    The mechanism is illustrated in \cref{fig:windows-merged}.

    To implement the above procedure,
    we present an ``outer'' algorithm \texttt{qcomp} in \cref{alg:qcomp} that ensures exact computation up to the
    specified precision $\out{w}$.
    The actual exact computation is specified by a pair of callback functions that are passed to \texttt{qcomp}.
    One of those functions sets up the precision and exponent, and the other performs the exact integer arithmetic
    for one element of the result vector.
    We need \gls*{bfp}-versions of the standard \gls*{blas} routines \texttt{axpby} ($z \gets \alpha x + \beta y$) and \texttt{gemv} ($z \gets \alpha Ax + \beta y$).
    For efficiency and simplicity, we use custom routines for the special cases $z \gets x - y$ (\texttt{sub})
    and $z \gets Ax$ (\texttt{spmv}).
    The callbacks used for the exact matrix-vector multiplication \texttt{spmv} and the corresponding algorithms for \texttt{axpby} and \texttt{gemv} are listed in
    \cref{sec:appendix-bfp-blas-algos}.
    We denote that a routine is wrapped by \texttt{qcomp} by prefixing it with $\texttt{q}$; for instance, we have
    $\texttt{qspmv}(\dots) = \texttt{qcomp}( \texttt{espmv-setup}, \texttt{espmv-row}, \dots )$.

    %%%%%%%%%%%%%%%%%%%%%%%%%%%%%%%%%%%%%%%%%%
    %%%%%% ---- QCOMP (outer loop) ---- %%%%%%
    %%%%%%%%%%%%%%%%%%%%%%%%%%%%%%%%%%%%%%%%%%

    %\begin{figure}
    %    \input{figures/alg-trunc}
    %\end{figure}
%    \begin{figure}
        \noindent
\begin{algorithm}[t]
    \begin{center}
        \begin{minipage}{\linewidth}
            \footnotesize
            \caption{$\texttt{Quantized BFP vector computation: } \texttt{qcomp}$}
            \label{alg:qcomp}
            \begin{algorithmic}[1]

                % Input arguments
                \Require $
                \texttt{setup-func},\
                \texttt{comp-func},\
                \texttt{input},\
                \out{w},\
                \tmp{w}: \unsignedpositive,\
                \gamma: \bfp^1,\
                n: \unsignedpositive$

                % \Ensure $0 < \out{w} \leq \tmp{w} \leq \bfpq{\exact{z}}, \ \gamma > 0,\ \gamma \text{ normalized}$
                \Ensure $0 < \out{w} \leq \tmp{w}, \ \gamma > 0,\ \gamma \text{ normalized}$

                \Variables $
                \texttt{setup},\
                \exact{z}: \bfp^1\ (\text{exact result}), \
                \tmp{z}:\bfp^n\ (\text{temporary result}), \newline
                \out{z}: \bfp^n\ (\text{normalized result}), \
                \tmp{\mu},\
                \tmp{\lambda},\
                \exact{\mu},\
                \lambda_{\textrm{out}} : \signed,\
                \texttt{overflow},\
                \texttt{underflow}: \texttt{bool}
                $

                \LineComment{Allocation and setup.}

                % \State $\exact{z},\ (\texttt{output-setup}) \gets \texttt{setup-func}(\texttt{input})$
                \State $\left(\exact{z},\ \texttt{setup} \right)  \gets \texttt{setup-func}(\texttt{input})$
                \Comment{Setup necessary variables for exact computations.}

                \State $\tmp{\mu} \gets \texttt{msb}(\bfpm{\gamma}) + \bfpe{\gamma} - \bfpe{\exact{z}}$
                \Comment{MSB of the estimated mantissa window.}

                \State $\tmp{\lambda} \gets \tmp{\mu} - \tmp{w}$
                \Comment{LSB of the estimated mantissa window.}

                \State $\bfpq{(\out{z})} \gets \out{w}$
                \Comment{Width of result BFP vector.}

                \State $\bfpq{(\tmp{z})} \gets \tmp{w}$
                \Comment{Width of temporary BFP vector.}

                \State $\bfpe{(\tmp{z})} \gets \bfpe{\exact{z}} + \tmp{\lambda}$
                \Comment{Block-exponent of truncated temporary result $\tmp{z}$.}

                \LineComment{Computation of temporary BFP vector.}

                \State $\exact{\mu} \gets 1$
                \Comment{Keep track of the largest MSB idx. of the exact result.}

                \ParFor{$i = 1,\, \dots,\, n$}

                \State $\exact{z} \gets \texttt{comp-func}(\exact{z},\ \texttt{input},\ \texttt{setup},\ i)$
                \label{alg:qcomp:line:firstloopexact}

                \Comment{Compute exact quantity in parallel.}

                \State $\exact{\mu} \gets \texttt{atomic-max}(\exact{\mu},\ \texttt{msb}(\bfpm{\exact{z}}))$
                \label{alg:qcomp:line:atomic}
                \Comment{Atomic update of the global max.}

                \State $\bfpm{\exact{z}} \gets \bfpm{\exact{z}} \RightShift \tmp{\lambda}$
                \Comment{Shift and truncate exact result.}

                \State $(\bfpm{(\tmp{z})})_i \gets \text{decr} \left( \bfpq{\exact{z}} - \bfpq{(\tmp{z})},\ \bfpm{\exact{z}} \right)$
                \label{alg:qcomp:line:firstlooptrunc}

                \Comment{Cast result into temporary BFP vector.}

                \EndParFor
                \label{alg:qcomp:line:endfirstloop}

                \LineComment{Ensuring exact result up to $\out{w}$ bits.}

                %\State $\exact{\mu} \gets \texttt{reduce-max}(\gamma_\text{max})$
                %\Comment{Find maximum over all processes.}

%                \If{$(\Gamma_\text{max} = 0)$}
%                \nk{fix this if stmt}
%                    \State\Return $(0, \mathbf{0}, \out{w}, n)$
%                    \Comment{Return a zero BFP vector.}
%                \EndIf

                %\State $s \gets \lfloor\log_2(R_\text{max})\rfloor - \bfpe{\exact{r}} \geq 0$
                %\Comment{Actual bit idx of MSB.}

                % \State $\exact{\mu} \gets \min( \lfloor \log_2(\Gamma_\text{max}) \rfloor - \bfpe{\exact{z}} +
                % 2,\ \bfpq{\exact{z}} + 1 )$
                %\State $\exact{\mu} \gets \Gamma_\text{max}$
                %\Comment{MSB of the actual mantissa.}

                \State $\out{\lambda} \gets \exact{\mu} - \out{w}$
                \Comment{LSB of the target mantissa window.}

                \State $\bfpe{(\out{z})} \gets \bfpe{\exact{z}} + \out{\lambda}$
                \Comment{Block-exponent of result $\out{z}$.}

                \State $\texttt{overflow} \gets {\tmp{\mu}} < {\exact{\mu}}$
                \Comment{Max. MSB index of exact result left of est. window.}

                \State $\texttt{underflow} \gets {\out{\lambda}} < {\tmp{\lambda}}$
                \Comment{Captured less than $\out{w}$ relevant bits.}

                % \State $\out{w} \gets \exact{\mu} - \out{\lambda}$
                % \Comment{Depending on $\exact{\mu}$ there may be less bits available than intended.}

                \If{$\texttt{overflow} \lor \texttt{underflow}$}\label{alg:qcomp-recomp}

                \LineComment{Overflow or less than $\out{w}$ meaningful bits. Recompute.}

                \ParFor{$i = 1,\, \dots,\, n$}

                % \State $\exact{z} \gets \texttt{comp-func}(\exact{z},\ \texttt{input},\
                % (\texttt{output-setup}),\ i)$
                \State $\exact{z} \gets \texttt{comp-func}(\exact{z},\ \texttt{input},\ \texttt{setup},\ i)$
                \label{alg:qcomp:line:recomputeexact}

                \Comment{Recompute exact quantity in parallel.}

                % \State $(\bfpm{(\out{z})})_i \gets \texttt{trunc}(\exact{z},\ \out{z},\ \lambda_\textrm{out},
                % \ i)$
                \State $(\bfpm{(\out{z})})_i \gets \text{decr} \left(\bfpq{\exact{z}} - \bfpq{(\out{z})},\ \bfpm{\exact{z}} \RightShift \out{\lambda} \right)$
                \label{alg:qcomp:line:recomputetrunc}

                \Comment{Shift and truncate exact result.}

                \EndParFor

                \Else

                    \LineComment{At least $\out{w}$ meaningful bits}

                    \ParFor{$i = 1,\, \dots,\, n$}

                    \State $(\bfpm{(\out{z})})_i \gets \text{decr}\left(\tmp{w} - \out{w},\ (\bfpm{(\tmp{z})})_i
                    \RightShift (\out{\lambda} - \tmp{\lambda})\right)$
                    \label{alg:qcomp:line:successtrunc}

                    \Comment{Truncate all bits that exceed the target window.}

                    \EndParFor

                \EndIf

                \State\Return $\out{z}$

            \end{algorithmic}
        \end{minipage}
    \end{center}
\end{algorithm}

%    \end{figure}

%%%%%%%%%%%%%%%%%%%%%%%%%%%%%%%%%%%%%%%%%%%%%%%%%%%%%%%%%%%%%%%%%%%%%%%%%%%%%%%%%%%%%%%%%%%%%%%%%%%%%%%%%%%%%%%%%%%%%%%%
%%%%%%%%%%%%%%%%%%%%%%%%%%%%%%%%%%%%%%%%%%%%%%%%%%%%%%%%%%%%%%%%%%%%%%%%%%%%%%%%%%%%%%%%%%%%%%%%%%%%%%%%%%%%%%%%%%%%%%%%
%%%%%%%%%%%%%%%%%%%%%%%%%%%%%%%%%%%%%%%%%%%%%%%%%%%%%%%%%%%%%%%%%%%%%%%%%%%%%%%%%%%%%%%%%%%%%%%%%%%%%%%%%%%%%%%%%%%%%%%%

    \section{Mixed- and progressive-precision multigrid}\label{sec:mg}

    As in ~\cite{Tamstorf:2021:DiscretizationErrorAccurateMixedPrecisionMultigrid},
    our goal is to approximate the solution of linear elliptic \glspl*{pde} up to discretization-error-accuracy using
    arithmetic of minimal precision.
    Thus, we consider linear systems of the form $Ax=b$
    with $A \in \real^{n \times n}$ SPD, $x, b \in \real^{n}$.
    A balance of quantization, discretization, and algebraic errors must be obtained by an appropriate choice of the
    precisions employed during computation.
    The system is solved by iterative refinement (\ir) with an inner V-cycle (\V), possibly as part of full
    multigrid (\fmg).
    We apply subscripts to relate a quantity to a refinement level (for example, $A_j$ refers to the discrete operator
    on level $j$).
    As in~\cite{Tamstorf:2021:DiscretizationErrorAccurateMixedPrecisionMultigrid},
    three precisions $\eweq_j \leq \ewe_j
    \leq \ewed_j$ are defined on each refinement level $j > 0$.
    The ``working'' precision, i.e.,~the precision of the computed result, is $\ewe$.
    It is used in \ir~and \fmg, while the precision of the inner solver is reduced to $\ewed$.
    To account for quantization errors induced by storing the input in finite precision, $A$ and $b$ are stored in
    $\eweq$ precision.
    A fourth, high precision $\eweb_j \leq \ewe_j$ is required
    in~\cite{Tamstorf:2021:DiscretizationErrorAccurateMixedPrecisionMultigrid}
    to ensure more precise computation of the residual in \ir.
    However, thanks to \cref{alg:qcomp}, we can assert that each result is computed exactly up to a specified precision,
    thereby allowing us to eliminate the need for $\eweb_j$.
    Under certain assumptions, \cite{Tamstorf:2021:DiscretizationErrorAccurateMixedPrecisionMultigrid} shows that the precisions required to attain discretization-error-accuracy can be
    bounded by functions of the following quantities: the finite element polynomial of order $k$, the order $2m$ of the \gls*{pde}, and the pseudo
    mesh size $h_j = \underline{\kappa}_j^{-\frac{1}{2m}}$, where
    $\underline{\kappa}_j \coloneqq \| |A_j| \| \cdot \|A_j^{-1}\|$ ($|\cdot|$ denotes matrix entries replaced by their
    absolutes).
    In particular, it is shown that
    \begin{equation}
        \label{eq:precision-bounds}
        \eweq_j \in \mathcal{O}\left(h_j^{k+m}\right), \quad
        \ewe_j \in \mathcal{O}\left(h_j^{k}\right), \quad
        \ewed_j \in \mathcal{O}\left(h_j^{m}\right).
    \end{equation}

    Following \Cref{sec:bfp} for the relation of floating point to \gls*{bfp} precision, we are now
    interested in the behavior of the total error with respect to the corresponding \gls*{bfp} mantissa widths
    $\weweq_j$, $\wewe_j$, and $\wewed_j$.
    Assuming that $2 h_{j+1} = h_j$ and that \cref{eq:precision-bounds} also applies to \gls*{bfp} arithmetic, then
    \cref{eq:precision-bounds} suggests that the widths of the corresponding mantissas are
    related to refinement by
    \begin{equation}
        \label{eq:width-bounds}
        \weweq_j \in \mathcal{O}\left((k+m)j\right), \quad
        \wewe_j \in \mathcal{O}\left(kj\right), \quad
        \wewed_j \in \mathcal{O}\left(mj\right).
    \end{equation}
    As an example, consider the solution of a second-order \gls*{pde} ($m = 1$).
    According to \cref{eq:width-bounds}, only $1$ bit needs to be added per refinement level to
    the mantissa width $\wewed_j$ used in the inner solver to ensure discretization-error-accuracy.

    \Cref{alg:ir,alg:v-cycle,alg:fmg} list the \gls*{bfp}-versions of \ir, \V, and \fmg~as defined in~\cite{Tamstorf:2021:DiscretizationErrorAccurateMixedPrecisionMultigrid}.
    We refer to \ir~with \V~as the inner solver by \ir-\V.
    The \gls*{bfp}-routines require the width of the mantissa for the result $\out{w}$, an estimate for the infinity norm
    of the result $\gamma$, and the number of bits $\tmp{w} \geq \out{w}$ to be used for the mantissa of the
    temporary result.
    Note that line~\ref{alg:ir-correction} in \cref{alg:ir} (correction step) and line~\ref{fmg-interp} in \Cref{alg:fmg} (\fmg-prolongation) are the only calls where we use standard working precision $\out{w} = \wewe$.
    For all remaining calls, we use low precision $\out{w} = \wewed$.
    We comment on the choice of $\gamma$ and $\tmp{w}$ in \Cref{sec:choice-gamma-wtmp,sec:skip-normalization}
    but omit them as input arguments
    to the \gls*{bfp}-routines
    in \cref{alg:ir,alg:v-cycle,alg:fmg} for better readability.

    Inside of \V~(\cref{alg:v-cycle}) we use a second-order Chebyshev iteration for relaxation. 
    We base our implementation on a simplification of~\cite[Algorithm 1]{Gutknecht2002}.
    When reduced to just two iterations and a zero initial guess, it can be implemented using a single call to \texttt{qgemv} as shown in line~\ref{alg:v-cycle-relaxation} of \cref{alg:v-cycle}.
    This call requires two coefficients $c_1, c_2$, which are estimated on refinement level $\ell_{\mathrm{est}} = 5$,
    using \cref{alg:cheby-nodes}.
    As in~\cite{Tamstorf:2021:DiscretizationErrorAccurateMixedPrecisionMultigrid}, the spectral radius $\rho(D^{-1}A)$
    is estimated via the solution of the generalized eigenvalue problem $A_{\ell_{\mathrm{est}}}x = \lambda D_{\ell_{\mathrm{est}}}x$,
    and the targeted percentage of the spectrum $\eta$ is determined empirically, by minimization of the V-cycle convergence rate
    over a set of values $\mu_i^* = i/100,\ i = 0, 1, \dots, 100$.
    After the computation of $c_1$ and $c_2$, the simplification also requires a setup phase where we set
    $A_j \gets D_j^{-1}A_j$, $b_j \gets D_j^{-1}b_j$, and $R_j \gets D_{j-i}^{-1} P^T D_{j}$ in order to avoid
    division operations (which are delicate in \gls*{bfp}-arithmetic)
    during the relaxation step. See \cref{alg:solver-setup}.
    The setup computations are assumed to be executed in exact arithmetic. (Concrete implementation details are
    described in \Cref{sec:results}.)

%%%%%%%%%%%%%%%%%%%%%%%%%%%%%%%%%%%%%%%%%%%%%%%%%%%%%%%%%%%%%%%%%%%%%%%%%%%%%%%%%%%%%%%%%%%%%%%%%%%%%%%%%%%%%%%%%%%%%%%%
%%%%%%%%%%%%%%%%%%%%%%%%%%%%%%%%%%%%%%%%%%%%%%%%%%%%%%%%%%%%%%%%%%%%%%%%%%%%%%%%%%%%%%%%%%%%%%%%%%%%%%%%%%%%%%%%%%%%%%%%
%%%%%%%%%%%%%%%%%%%%%%%%%%%%%%%%%%%%%%%%%%%%%%%%%%%%%%%%%%%%%%%%%%%%%%%%%%%%%%%%%%%%%%%%%%%%%%%%%%%%%%%%%%%%%%%%%%%%%%%%

    \section{Numerical results}\label{sec:results}

    Using the ideas described above, 
    \Cref{sec:results-minimal-precision} presents numerical results for two model problems that suggest that the precision bounds in
    \cref{eq:precision-bounds} also apply in \gls*{bfp}-arithmetic.
    \Cref{sec:results-estimated-precision} covers the a priori estimation of the individual \gls*{bfp} precisions $\weweq_j, \wewe_j, \wewed_j$,
    and compares multigrid convergence rates of fixed-precision floating point, fixed-precision \gls*{bfp}, and progressive-precision \gls*{bfp} setups.
    We also report on the estimation of $\gamma$ and $\tmp{w}$ to avoid recomputations in \cref{alg:qcomp} (\Cref{sec:choice-gamma-wtmp}),
    and study the effect of skipping \gls*{bfp}-vector normalization altogether (\Cref{sec:skip-normalization}).

    With $\Omega \coloneqq (0, 1)$ and $f \in L^2\left( \Omega^d \right)$, we consider the following model problems:
    find $u \in C^{2m}$ s.t.

    \noindent\begin{minipage}{0.49\textwidth}
                 \begin{equation}
                     \label{eq:pde-poisson}
                     \begin{aligned}
                         - \Delta u &= f \quad \text{in $\Omega^d$}, \\
                         u &= 0 \quad \text{on $\partial\left(\Omega^d\right)$},
                     \end{aligned}
                 \end{equation}
    \end{minipage}
    \begin{minipage}{0.49\textwidth}
        \begin{equation}
            \label{eq:pde-biharmonic}
            \begin{aligned}
                u'''' &= f \quad \text{in $\Omega$}, \\
                u = u' &= 0 \quad \text{on $\partial\Omega$},
            \end{aligned}
        \end{equation}
    \end{minipage}
    \medskip

    \noindent where $d \in \{1, 2\}$, $m = 1$ in \cref{eq:pde-poisson}, and $d = 1$, $m = 2$ in \cref{eq:pde-biharmonic}.
    The biharmonic equation $\cref{eq:pde-biharmonic}$ is selected due to the rapidly growing condition number of the
    system matrix of the discrete problem, which is especially challenging for low-precision
    computations~\cite{Tamstorf:2021:DiscretizationErrorAccurateMixedPrecisionMultigrid}.
    Both model problems are approximated via the standard Rayleigh-Ritz finite element method, using identical,
    finite-dimensional trial and test spaces.
    For the discetization, we use B-spline finite elements of order $k = p + 1$, where $p$ is the polynomial degree.
    The Dirichlet boundary conditions are enforced strongly.
    Overall, the setup follows~\cite{Tamstorf:2021:DiscretizationErrorAccurateMixedPrecisionMultigrid}.
    The manufactured solutions $u$ are chosen as smooth functions with trigonometric components.
    All setup computations (including assembly of the linear system and \cref{alg:solver-setup}) are performed
    in high precision floating point arithmetic using a $400$ bit mantissa to ensure sufficient accuracy.
    For comparison, double precision has a 53 bit mantissa and quad precision a 113 bit mantissa.
    Integrals are approximated via Gauss-Legendre quadrature, with $(p + 1)^{d}$ nodes and weights per element.

    For all experiments, we use a prototype C++ \gls*{bfp} implementation based on the GNU Multiple Precision
    Arithmetic Library (GMP) offering arbitrarily wide integer types, and the GNU Multiple Precision Floating-Point
    Reliable Library (MPFR) for arbitrary precision floating point formats.
    Our implementation is experimental and favors flexibility over computational performance as the numerical results
    are the focus of this paper.
    Thus, we do not present any run time results.

    \subsection{Confirmation of BFP-precision bounds}\label{sec:results-minimal-precision}

    It is desirable to estimate sufficient mantissa widths $\weweq_j$, $\wewe_j$, and $\wewed_j$ before application
    of the solver.
    Given a specific problem, the asymptotic bounds in \cref{eq:precision-bounds,eq:width-bounds} are used for such
    a priori estimates.
    The objective of this section is the experimental confirmation of those bounds, to assert that they can in fact
    be used for a priori estimates in practice.
    To that end, the mantissa widths $\weweq_j$, $\wewe_j$, and $\wewed_j$ are initially not estimated, but iteratively increased in
    steps of one bit per run of \ir-\V~(starting from 1 bit) until the \gls*{bfp}-approximation is close to a reference solution for each level $j$.
    We compare the computed (\gls*{bfp}-) solution $\tilde{u}_h$ to a reference $u_h$ that is computed in floating
    point arithmetic using a $400$ bit mantissa.
    We accept $\tilde{u}_h$ if $\|u - \tilde{u}_h\|_\mathcal{L} / \|u - u_h\|_\mathcal{L} \leq 1.5$, where
    $\|\cdot\|_\mathcal{L} = a( \cdot, \cdot )^{1/2}$ is the energy norm and $a$ is the bilinear form associated with the
    weak formulation.
    We aim for discretization accuracy, i.e.,~a total error of order $\mathcal{O}(h^{k-m})$.
    This is achieved for the reference solution.
    The prolongated exact solution of the next coarser grid is used as an initial guess, mimicking \fmg, and
    the number of \ir-iterations is limited to 50.
    This is an overly pessimistic limit for \fmg~and, in most cases, a few iterations are sufficient.

    In this initial experiment, we are interested in finding the smallest mantissa required to observe convergence and not in the
    convergence rate of the multigrid solver.
    Precisions that are sufficient to achieve discretization-error-accuracy
    do not necessarily lead to satisfactory convergence rates of the linear solver.
    Increasing the precision beyond what is required to achieve discretization-error-accuracy may further increase
    the convergence rates.
    Clearly, this rate is limited, and higher precision generally entails lower computational performance.
    In practice, a trade-off has to be made, which we revisit in \Cref{sec:results-estimated-precision}.

\begin{figure}
        \centering
        \resizebox*{\textwidth}{!}{%
            \input{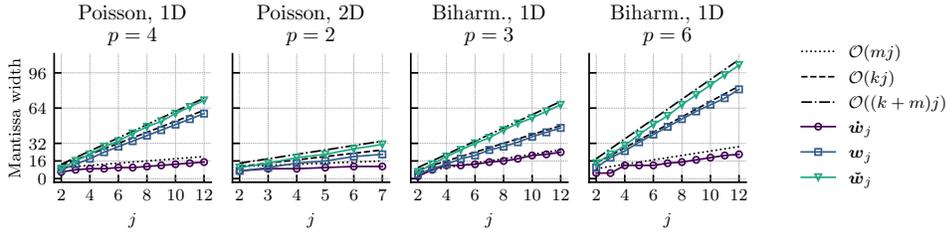}%
        }
        \vspace{-24pt}
        \caption{Minimum number of bits required for the \gls*{bfp} mantissa
        to maintain optimal error convergence in the energy norm over grid refinement for four test cases using \ir-\V.
        The x-axis shows the refinement level with mesh sizes $h = 2^{-j}$.}
        \label{fig:bfp-min-width}
    \end{figure}

    For simplicity in this section, we do not employ progressive precision inside the V-cycle itself, but apply the precision of
    the finest grid throughout the hierarchy for each run.
    First, we choose $\wewe_j = \wewed_j = 200$ and determine the minimal width
    $\weweq_j$ over mesh refinement such that $\tilde{u}_h$ fulfills the convergence criterion.
    Using the obtained precisions $\weweq_j$, we fix $\wewed_j = 200$ and find the
    minimal $\wewe_j$ in the same way.
    Eventually using both $\weweq_j$ and $\wewe_j$, we apply the same process to find
    the minimal $\wewed_j$.
    The results for four test cases are plotted in \cref{fig:bfp-min-width}.

    We observe that the asymptotic behavior in \cref{eq:width-bounds} that is predicted for floating point arithmetic
    in~\cite{McCormick:2021:AlgebraicErrorAnalysis,Tamstorf:2021:DiscretizationErrorAccurateMixedPrecisionMultigrid}
    also holds for the \gls*{bfp} implementation.
    Note that the finite element polynomial order $k$ is $p+1$, where $p$ is the polynomial degree.
    For the biharmonic equation with $p=3$, we have $m = 2$, $k = 4$, and therefore expect
    (and observe)
    $\weweq_j \in \mathcal{O}(6j)$,
    $\wewe_j \in \mathcal{O}(4j)$, and
    $\wewed_j \in \mathcal{O}(2j)$.
    We do not see any asymptotic difference between the 1D and 2D cases (nor do we expect any), but we include the 2D
    Poisson test case to illustrate this. 

    \subsection{A priori estimation of required BFP-precision}\label{sec:results-estimated-precision}

    Selecting the precisions via the approach in \Cref{sec:results-minimal-precision} is expensive and impractical for
    real applications.
    % Therefore, we employ the precision estimation algorithm
    % from~\cite{Tamstorf:2021:DiscretizationErrorAccurateMixedPrecisionMultigrid}.
    It is feasible, however, to estimate the precisions based on discretization-dependent constants
    that can be computed relatively cheaply on very coarse grids.
    Those estimated precisions are then extrapolated for finer grids according to their asymptotic behavior given in
    \cref{eq:precision-bounds,eq:width-bounds}.

    In this paper, we choose the precisions as follows.
    For $\wewe_j$, we employ the estimation algorithm
    from~\cite{Tamstorf:2021:DiscretizationErrorAccurateMixedPrecisionMultigrid}.
    This algorithm could also be used with small adjustments to determine $\wewed_j$ and $\weweq_j$, but instead we
    propose to select $\wewed_j$ and $\weweq_j$ subject to a certain target convergence rate of the solver.
    % \rt{All of the following should be encoded in an algorithm.}
    For that, we employ \cref{alg:bfp-prec-est} with $j_c = 5$, $q_{\mathrm{max}} = 64$, and $\rho_{\mathrm{thresh}} = 1.05$.
    The minima are determined using binary search over $\{1, \dots, q_{\mathrm{max}}\}$.
    On coarse levels $j_c$, this approach is reasonably fast.
    \begin{algorithm}[h]
    \begin{center}
        \begin{minipage}{\linewidth}
            \footnotesize
            \caption{$\texttt{BFP precision estimation: } \texttt{bfp-prec-est}$}
            \label{alg:bfp-prec-est}
            \begin{algorithmic}[1]

                % Input arguments
                \Require
                $j_c : \unsignedpositive$ (estimation level),
                $q_{\mathrm{max}} : \unsignedpositive$ (sufficient mantissa width for convergence),
                $\rho_{v,\mathrm{thresh}} : \real$ (threshold for relative convergence rate),
                $m : \unsignedpositive$ ($2m = $ order of the \gls*{pde}),
                $k : \unsignedpositive$ (approximation order),
                $\wewe_j : \unsignedpositive,\, j = 1, \dots, \ell$ (precisions estimated as
                    in~\cite{Tamstorf:2021:DiscretizationErrorAccurateMixedPrecisionMultigrid})

                %\Ensure $.$

                %\Variables $.$

                \State $\weweq_j(q) \coloneqq j(m + k) + q$
                \Comment{Shorthand for \cref{eq:width-bounds} plus constant.}

                \State $\wewed_j(q) \coloneqq jm + q$
                \Comment{Shorthand for \cref{eq:width-bounds} plus constant.}

                \State $\rho_{v,\mathrm{ref}} \gets \texttt{conv-rate-v-cycle}(
                    \weweq_j(q_{\mathrm{max}}), \wewe_j, \wewed_j(q_{\mathrm{max}}), j_c
                )$

                \Comment{Reference convergence rate on level $j_c$.}

                \State $\check{q} \gets \min \left\{
                    q \in \{1, \dots, q_{\mathrm{max}}\} :
                    \texttt{conv-rate-v-cycle}(
                        \weweq_j(q), \wewe_j, \wewed_j(q_{\mathrm{max}}), j_c
                    ) / \rho_{v,\mathrm{ref}} < \rho_{v,\mathrm{thresh}}
                \right\}$

                \Comment{Min. additive constant for $\weweq_j$ to satisfy conv. crit.}

                \State $\dot{q} \gets \min \left\{
                    q \in \{1, \dots, q_{\mathrm{max}}\} :
                    \texttt{conv-rate-v-cycle}(
                        \weweq_j(\check{q}), \wewe_j, \wewed_j(q), j_c
                    ) / \rho_{v,\mathrm{ref}} < \rho_{v,\mathrm{thresh}}
                \right\}$

                \Comment{Min. additive constant for $\wewed_j$ to satisfy conv. crit.}

                \State\Return $(\weweq_j(\check{q}), \wewed_j(\dot{q}))$

            \end{algorithmic}
        \end{minipage}
    \end{center}
\end{algorithm}

    The convergence rate $\rho_{v} = \|V\|_A$ is computed as the square root of the largest generalized eigenvalue of
    $V^T A V x = \lambda A x$, where $V$ is the error propagation matrix of \V.
    This matrix is constructed by applying \V~to the canonical basis vectors.
    This is the same approach as taken in~\cite{Tamstorf:2021:DiscretizationErrorAccurateMixedPrecisionMultigrid}.
    In \cref{alg:bfp-prec-est}, $\texttt{conv-rate-v-cycle}(\weweq_j, \wewe_j, \wewed_j, j)$ computes $\rho_{v}$
    on level $j$, using progressive precision to construct $V$ in \gls*{bfp} arithmetic.

\begin{figure}
        \centering
        \resizebox*{\textwidth}{!}{%
            \input{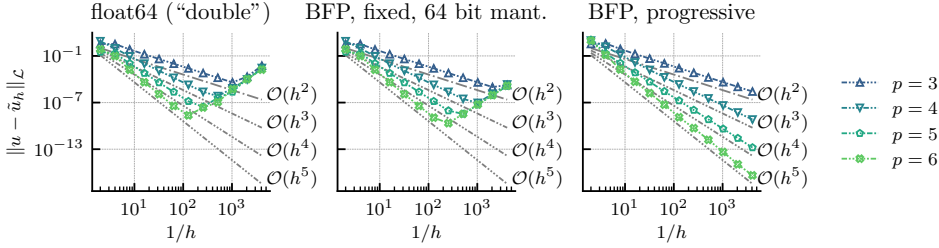}
        }
        \vspace{-24pt}
        \caption{Total error over refinement applying \fmg~to the biharmonic model problem using different setups.
        The left plot shows the error evolution using standard 64 bit floating point precision; the plot in the
        center shows the results for \gls*{bfp}-\fmg, with fixed precision, i.e., $\weweq_j = \wewe_j = \wewed_j = 64$
        for all levels $j = 1, \dots, 12$. Both fixed precision approaches lead to approximations that are eventually
        dominated by rounding errors. The plot on the right employs the progressive precision \gls*{bfp}-\fmg~algorithm
        with estimated precisions. The corresponding computed solutions are discretization-error-accurate, regardless
        of the refinement level. As shown in~\cite{Tamstorf:2021:DiscretizationErrorAccurateMixedPrecisionMultigrid},
        the same qualitative results are also achieved for progressive precision floating point implementations.}
        \label{fig:bfp-fmg-vs-reference}
    \end{figure}

\begin{figure}
        \centering
        \resizebox*{\textwidth}{!}{%
            \input{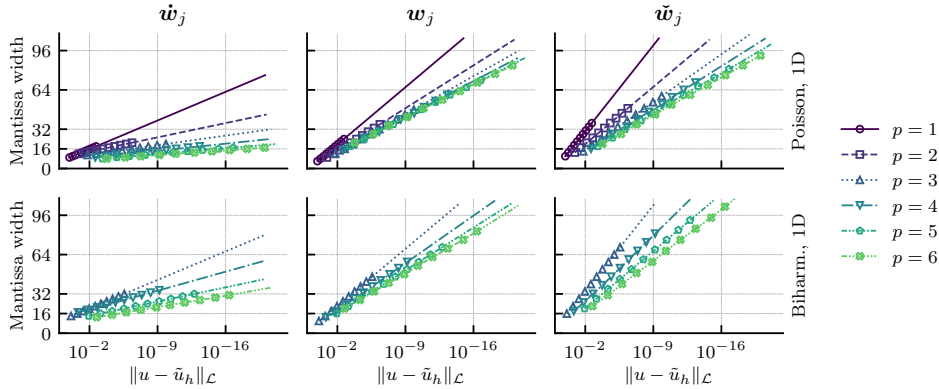}
        }
        \vspace{-24pt}
        \caption{Progressive precision \gls*{bfp}-\fmg~using a priori estimated mantissa widths.
        The plots show the applied precisions to achieve a certain target accuracy.
        Markers indicate actual data points.
        All lines have been extrapolated.
        The number of \ir-\V-iterations, $N$, per \fmg-level has been estimated according
        to~\cite{Tamstorf:2021:DiscretizationErrorAccurateMixedPrecisionMultigrid}
        as $N = \{2, 1, 1, 3, 7, 15\}$ for the Poisson problem with $p = \{1, \dots, 6\}$ and as
        $N = \{2, 1, 2, 4\}$ for the biharmonic equation with $p = \{3, \dots, 6\}$.
        The estimates for $N$ match the estimates that we obtained from a reference implementation
        that employs extremely accurate floating point arithmetic with mantissas using 400 bits. 
        (IEEE-754 ``double'' precision uses 53 bit mantissas.) The number of iterations for higher polynomial
        degrees are larger than necessary in practice.
        For each data point in x-direction, a refinement step has been performed.
        The rightmost data points for each line correspond to a mesh size of
            $h \approx \num{2.5e-4}$.}
        \label{fig:bfp-fmg}
    \end{figure}

    In \cref{fig:bfp-fmg-vs-reference}, the progressive precision \gls*{bfp}-\fmg~solver (right plot) is applied to the biharmonic
    equation and compared to a reference implementation using standard ``double'' precision 64 bit IEEE-754
    floating point arithmetic (left plot), and a \gls*{bfp}-\fmg~solver with fixed precision (center plot) on all levels.
    The results demonstrate the necessity of progressive precision, for both
    standard floating point and \gls*{bfp} implementations.
    For both fixed-precision \fmg~solvers, we applied $N = 20$ \ir-\V~iterations per level, which should be more than
    sufficient given the much smaller number of iterations that we need to achieve discretization-error-accuracy
    with the progressive precision \gls*{bfp} implementation. (The number of iterations for progressive precision
    \gls*{bfp}-\fmg~is listed in \cref{fig:bfp-fmg}.)
    Note that the fixed precision \gls*{bfp} setup even performs slightly better than the floating point version for this
    test case.
    This can be explained by the different usage of the 64 bits: the 64 bit floating point format reserves only 53 bits
    for the mantissa, while the \gls*{bfp} format uses all 64 bits for the mantissa.

    \Cref{fig:bfp-fmg} shows the mantissa widths used to achieve a certain target accuracy in the energy norm when using
    progressive precision \gls*{bfp}-\fmg.
    Asymptotically optimal grid convergence is observed under refinement for all cases.
    For the biharmonic equation, this is depicted in the right plot in \cref{fig:bfp-fmg-vs-reference}, which uses the same data.
    A comparison to the results of~\cite[Figures 2, 9]{Tamstorf:2021:DiscretizationErrorAccurateMixedPrecisionMultigrid}
    shows that the number of bits required to attain a certain target accuracy is similar to that of floating point
    arithmetic.

    \subsection{Choice of \texorpdfstring{$\gamma$}{gamma} and \texorpdfstring{$\tmp{w}$}{w}}\label{sec:choice-gamma-wtmp}

    \Cref{alg:qcomp} (\texttt{qcomp}) requires the estimation of an upper bound $\gamma$ of the infinity norm
    $\|z\|_\infty$ for the result $z$ of \cref{eq:gemv}, and the estimation of a sufficiently large mantissa width, $\tmp{w}$, for the
    temporary result.
    The objective is to choose these values as tight as possible, while still avoiding underflow and overflow along with 
    the ensuing recomputation (see \cref{alg:qcomp}, line~\ref{alg:qcomp-recomp}).

    In the following we discuss heuristics for choosing $\gamma$ and $\tmp{w}$ for the individual steps of
    \cref{alg:ir,alg:v-cycle,alg:fmg}, with a focus on estimates for \gls*{bfp}-\fmg.
    However, since $\tmp{w} \geq \out{w}$ we provide $w^*_{\mathrm{add}} = \tmp{w} - \out{w}$ where $\out{w}$ is level dependent while $w^*_{\mathrm{add}}$ is fixed across all levels.

    To derive the heuristics, we make several simplifying assumptions that characterize an important class of problems, but hopefully carry over to more general cases. In particular, assume in the following that $\diag(A)=I$ and that the entries of $E \coloneqq I-A$ and the row sums of $A$ are nonnegative. Note then that $\|E\|_\infty \leq 1$ and $\rho(A) \leq 2$.
    % \sm{Nils, please check to see if we really need this: Assume for further simplicity that $\rho(A)$ is actually equal to $2$.}
    The \ir-\V-iteration on one level of \fmg~is indicated with superscript $^{(i)},\ i \in \{1, \dots, N\}$, if necessary.
    Applying \ir-\V~directly (without \fmg) may involve slightly different assumptions, but the overall approach is similar
    (see \Cref{sec:skip-normalization,rem:skipping-normalization-ir-v}).

    The values of $\gamma$ are based on various bounds and inequalities, while the values of $w^*_{\mathrm{add}}$ are determined empirically. In particular, $w^*_{\mathrm{add}}$ is chosen such that no recomputation is
    triggered during \gls*{bfp}-\fmg~on the finest level $j = 12$ for both 1D model problems and $p \in \{1, \dots, 6\}$
    using the setup discussed in \Cref{sec:results-estimated-precision}. The same value of $w^*_{\mathrm{add}}$ is then later used for all levels.
    The choices for $\gamma$ and $w^*_{\mathrm{add}}$ are summarized in \cref{tab:estimates}.
    \paragraph{\ir, residual (\cref{alg:ir}, line~\ref{alg:ir-residual})}
    For a convergent \ir~process we generally expect the residual to steadily get smaller over time, but there is no guarantee that this will happen monotonically (which can lead to overflow), and occasionally large reductions may occur (which can lead to underflow). One could try to account for all of this, but we have found it to be best simply to choose
    $\tmp{w}$ a little larger that $\out{w}$ and assume that the residual does not increase. 
    To be more specific, for iteration $i$, we choose
    \begin{equation}
            \gamma =
            \begin{cases}
                 \|r_{\ell-1}^{(N)}\|_\infty, & i = 1, \\
                 \|r_{\ell}^{(i-1)}\|_\infty, & i > 1, \\
            \end{cases} \qquad
            w^*_{\mathrm{add}} =
            \begin{cases}
                 5, & i = 1, \\
                 4, & i > 1, \\
            \end{cases}
    \end{equation}
    where $r_{\ell-1}^{(N)}$ is the residual after \ir-\V~iteration $N$ on the next coarser \fmg-level.
    \paragraph{\ir, correction (\cref{alg:ir}, line~\ref{alg:ir-correction})}
    The bound for $\gamma$ in this step follows from the triangle inequality. We choose $\gamma = \|x\|_\infty + \|y\|_\infty$,
    which works extremely well in practice, so that no bits have to be added, i.e., $w^*_{\mathrm{add}} = 0$.
    \paragraph{\V, relaxation (\cref{alg:v-cycle}, line~\ref{alg:v-cycle-relaxation})}
     Second-order Chebyshev relaxation has the form
    \begin{equation}
        y \gets (c_1 I + c_2 A) r,
    \label{2nd-cheby}
    \end{equation}
    where $c_1$ and $c_2$ are scalar constants. Noting in \cref{alg:cheby-nodes} that $\alpha > 0, c > 0$, and
    \[
    \beta = (1 + \frac{c}{\sqrt{2}\alpha})(\alpha - \frac{c}{\sqrt{2}}) > \alpha - \frac{c}{\sqrt{2}} > \alpha - c = \eta\rho > 0
    \]
    it follows that $c_1 > 0$ and $c_2 < 0$. 
    Since $0 \le A \le 2I$, then $(c_1 + 2c_2) I \le c_1 I + c_2 A < c_1 I$, which in turn implies that
    \begin{equation}
        \|y\|_\infty / \|r\|_\infty \in [c_1 + 2c_2, c_1).
    \end{equation}
    We thus set $\gamma = c_1 \|r\|_\infty$.
    Using the lower bound for $\|y\|_\infty / \|r\|_\infty$, $\tmp{w}$ could be chosen as
    \begin{equation}
        \tmp{w} = \out{w} + \left\lceil \log_2 \left( \frac{c_1}{c_1+2c_2} \right)\right\rceil.
    \end{equation}
    In practice, simply adding $w^*_{\mathrm{add}} = 2$ bits works well in our experience.
    \paragraph{\V, residual (\cref{alg:v-cycle}, line~\ref{alg:v-cycle-residual})}
    Rewriting (\ref{2nd-cheby}) as
    \[
        y \gets c_2 Ar + c_1 r = c_2 (A-I)r + (c_1 - c_2)r = -c_2 Er + (c_1 - c_2)r
    \]
    leads to the following bound on the subsequent relative residual norm:
    \begin{equation}
        \begin{aligned}
            \|Ay - r\|_\infty / \|r\|_\infty \gets& \|-c_2 A E r + (c_1 - c_2) Ar - r\|_\infty / \|r\|_\infty \\
            &\leq \|-c_2 AE + (c_1 - c_2)A - I\|_\infty \\
            &\leq 2c_2 + 2 |c_1 - c_2| + 1 \\
            &= 2c_1+1,
        \end{aligned}
    \end{equation}
    where the last line follows because $c_1>c_2$.
    This gives a liberal upper bound for the range, that is, a conservative estimate for $\gamma$.
    We therefore scale the bound empirically by $1/4$, and choose $\gamma = (1/4) (2c_1 + 1)\|r\|_\infty$.
    Unfortunately, there is no useful lower bound for $\|Ay - r\|_\infty/\|r\|_\infty$ because we cannot rule out the possibility that
    the error $e=y - A^{-1}r$ is very smooth (e.g., the minimal eigenvector), meaning that $\|Ay - r\|_\infty/\|r\|_\infty$ would be $\mathcal{O}(h^2)$.
    There seems to be little choice here but to use an initial $\tmp{w} \gg \out{w}$
    and adjust it to a more conservative value based on the observed $y$ as the cycles proceed.
    In practice, adding a few bits to $\tmp{w}$ compared to $\out{w}$, however, works well.
    The hope is that the size of the residual after relaxation is somewhat consistent from one cycle to the next.
    We find $w^*_{\mathrm{add}} = 4$ to be sufficient.
    \paragraph{\V, restriction (\cref{alg:v-cycle}, line~\ref{alg:v-cycle-restriction})}
    The upper bound for the residual transfer
    % is much clearer \rt{Than what ?}, and
    follows from the triangle inequality, i.e., $\gamma = \|R\|_\infty \|r_v\|_\infty$.
    Here too we lack a useful lower bound on the range.
    Indeed, $R r_v = 0$ is certainly possible.
    However, such a loss indicates that coarsening is really useless itself. Fortunately it is also unnecessary because the residual
    must be oscillatory and relaxation alone would have reduced it significantly.
    That is, it is probably sufficient to choose $\tmp{w}$ only a little larger than $\out{w}$ because little damage would
    be done by any loss of bits resulting from an inaccurate estimate.
    However, to avoid underflow in all cases, we find that we need to add $w^*_{\mathrm{add}} = 6$ bits.
    \paragraph{\V, interpolation and corrrection (\cref{alg:v-cycle}, line~\ref{alg:v-cycle-correction})}
    For standard nodal-based interpolation, we can expect $\|Pd_{\ell-1}\|_\infty = \|d_{\ell-1}\|_\infty$, which gives $\gamma$
    (and $\tmp{w}$) exactly.
    Other forms of interpolation are probably at least approximately the same.
    For the correction, $y \gets y - d$, an obvious upper bound on the range is again given by the triangle inequality.
    The range is less clear except for the initial correction when $y=0$ so that the updated $y$ has norm $\|y\|_\infty=\|d\|_\infty$.
    For later cycles, especially near convergence, a major reduction in the size of $y$ would not be expected.
    This heuristic turns out to apply, and we set $\gamma = \|y\|_\infty + \|d_{\ell-1}\|_\infty$, only adding $w^*_{\mathrm{add}} = 1$ bit.
    \paragraph{\fmg, interpolation (\cref{alg:fmg}, line~\ref{alg:fmg-interpolation})}
    Same argument as for interpolation in \V. We set $\gamma = \|x\|_\infty$ and require no additional bits ($w^*_{\mathrm{add}} = 0$).

    \begin{table}
        \centering
        \footnotesize
        \caption{Concrete estimates employed for $\gamma$, and experimentally determined values for
            $w^*_{\mathrm{add}}$ for \cref{alg:qcomp} (\texttt{qcomp}) in \gls*{bfp}-\fmg.
            The number of additional bits $w^*_{\mathrm{add}} = \tmp{w} - \out{w}$ is chosen minimally such that
            no recomputation is trigerred on the finest level $j = 12$ for both 1D model problems and $p \in \{1, \dots, 6\}$.
            The individual heuristics are discussed in \Cref{sec:choice-gamma-wtmp}.
            $N$ and $i \in \{1, \dots, N\}$ are used to indicate the number of \ir-iterations and the current
            \ir-iteration, if necessary.
        }
        \label{tab:estimates}
        \renewcommand{\arraystretch}{1.2}
        \begin{tabular}{lll}
            \toprule
            Step & $\gamma$ & $w^*_{\mathrm{add}}$ \\ % & \makecell{$\frac{w^*_{\mathrm{add}}}{\out{w}}$ (Bihm., $p = 6$, $j = 12$) }\\
            \midrule
            \ir, residual (line~\ref{alg:ir-residual})&
            $\begin{cases}
                 \|r_{\ell-1}^{(N)}\|_\infty, & i = 1 \\
                 \|r_{\ell}^{(i-1)}\|_\infty, & i > 1 \\
            \end{cases}$ &
            $\begin{cases}
                 5, & i = 1 \\
                 4, & i > 1 \\
            \end{cases}$
            \\
            \ir, correction (line~\ref{alg:ir-correction})&
            $\|x\|_\infty + \|y\|_\infty$ &
            $0$
            \\
            \V, relaxation (line~\ref{alg:v-cycle-relaxation})&
            $c_1 \|r\|_\infty$ &
            $2$
            \\
            \V, residual (line~\ref{alg:v-cycle-residual})&
            $\frac{1}{4} (2c_1 + 1) \|r\|_\infty$ &
            $4$
            \\
            \V, restriction (line~\ref{alg:v-cycle-restriction})&
            $\|R\|_\infty \|r_v\|_\infty$ &
            $6$
            \\
            \V, correction (line~\ref{alg:v-cycle-correction}) &
            $\|y\|_\infty + \|d_{\ell - 1}\|_\infty$ &
            $1$
            \\
            \fmg, prolongation (line~\ref{fmg-interp})&
            $\|x\|_\infty$ &
            $0$
            \\
            \bottomrule
        \end{tabular}
    \end{table}

    \begin{table}
        \centering
        \footnotesize
        \caption{Number of calls to $\texttt{qcomp}$ (\cref{alg:qcomp}) that triggered recomputation during \gls*{bfp}-\fmg~on the
        finest level $j = 12$, using $\tmp{w} = \out{w} + \min\left( w^*_{\mathrm{add}}, w^{\mathrm{max}}_{\mathrm{add}} \right)$.
        The values of $w^*_{\mathrm{add}}$ are listed in \cref{tab:estimates} and are chosen minimally such that no
        recomputation is necessary.
        The total number of calls to $\texttt{qcomp}$ on each level $j > 1$ is shown in parentheses.}
        \label{tab:recomputation}
        \begin{tabular}{llrrrrrrr}
            \toprule
            PDE & $w^{\mathrm{max}}_{\mathrm{add}}$ & $p = 1$ & $p = 2$ & $p = 3$ & $p = 4$ & $p = 5$ & $p = 6$ \\
            \midrule
            Poisson 1D & $\infty$ & $0$ $(13)$ & $0$ $(7)$ & $0$ $(7)$ & $0$ $(19)$ & $0$ $(43)$ & $0$ $(91)$ \\
             & $4$ & $0$ $(13)$ & $0$ $(7)$ & $0$ $(7)$ & $0$ $(19)$ & $0$ $(43)$ & $0$ $(91)$ \\
             & $2$ & $2$ $(13)$ & $2$ $(7)$ & $2$ $(7)$ & $3$ $(19)$ & $1$ $(43)$ & $4$ $(91)$ \\
             & $0$ & $9$ $(13)$ & $3$ $(7)$ & $4$ $(7)$ & $11$ $(19)$ & $24$ $(43)$ & $47$ $(91)$ \\
            \addlinespace[6pt]
            Biharm. 1D & $\infty$ & - & - &  $0$ $(13)$ & $0$ $(7)$ & $0$ $(13)$ & $0$ $(25)$ \\
             & $4$ & - & - &  $0$ $(13)$ & $1$ $(7)$ & $0$ $(13)$ & $2$ $(25)$ \\
             & $2$ & - & - &  $4$ $(13)$ & $1$ $(7)$ & $4$ $(13)$ & $5$ $(25)$ \\
             & $0$ & - & - &  $8$ $(13)$ & $3$ $(7)$ & $8$ $(13)$ & $13$ $(25)$ \\
            \bottomrule
        \end{tabular}
    \end{table}

    One observation is that the heuristics for $\gamma$ perform quite well.
    Although only a few ($\leq 6$) bits are added for the temporary vector, recomputations can be avoided altogether.
    The correction in \ir~and the \fmg~prolongation remarkably do not require additional bits to prevent underflow.
    Also, for the \V-cycle correction and the relaxation, only a respective 1 and 2 additional bits are sufficient.
    The results of the residual computation in \ir~and inside the \V-cycle, as well as the restriction, are less predictable.
    However, even for those computations, adding $4$ bits still yields good overall results.

    The bottom line is that, given some relatively simple estimates for the upper bounds, by adding only a few bits
    compared to $\out{w}$, recomputations can largely be avoided in \gls*{bfp}-\fmg. Furthermore, the exact choice of how many bits to add is not critial. To illustrate this, \cref{tab:recomputation} shows the number of recomputations triggered for other, smaller choices of
    $w^*_{\mathrm{add}}$.
    Concretely, we show results for setting
    $\tmp{w} =  \out{w} + \min\left(w^*_{\mathrm{add}}, w^{\mathrm{max}}_{\mathrm{add}}\right)$
    for $w^{\mathrm{max}}_{\mathrm{add}} \in \{0, 2, 4\}$.
    In practice, simply choosing $w^*_{\mathrm{add}} = 2$ for all operations yields a fairly small number of recomputations.

    \subsection{Skipping \texorpdfstring{\gls*{bfp}}{BFP}-vector normalization}\label{sec:skip-normalization}

    \Cref{alg:qcomp} ensures that the results of \gls*{bfp} vector-vector and matrix-vector operations are exact up to
    the specified precision via normalization of the computed \gls*{bfp}-vector.
    The results above indicate that with proper estimates for $\gamma$ and $\tmp{w}$, it provides a framework to build
    efficient discretization-error-accurate multigrid solvers in \gls*{bfp} arithmetic.
    \Cref{tab:recomputation} suggests that the estimates are, in fact, relatively accurate, and thus
    lead to the question whether normalization is necessary in the first place.
    To explore this, we elaborate on the accuracy of $\gamma$.

    In the best case, that is, when the upper bound $\gamma$ to the infinity norm of the result vector is estimated correctly
    ($\gamma = \|z\|_\infty$),
    the second part of \cref{alg:qcomp}, after line~\ref{alg:qcomp:line:endfirstloop}, has no effect on the output.
    Consequently, just setting $\out{w} = \tmp{w}$ and returning after the first loop would yield the same result,
    but it avoids the second pass over the vector.
    So the global maximum \gls*{msb} index $\exact{\mu}$ is not required, and the atomic update in
    line~\ref{alg:qcomp:line:atomic} can be skipped as well.
    Depending on the block size, which for our studies is chosen maximally, avoiding the atomic update may have
    significant (positive) impact on (parallel) performance.

    If the estimate for $\gamma$ is a little too large, and we return early, some precision is lost, i.e., the
    result is computed accurately only to less than $\out{w}$ bits.
    \Cref{tab:recomputation} shows that this just truncates the rightmost
    $\tmp{w} - \out{w}$ bits, i.e., $6$ bits in extreme cases, but mostly much less than that.
    The hope is that this precision loss is insignificant and that it still leads to a discretization-error-accurate method,
    possibly at the cost of a slightly worse solver convergence rate.

    If the estimate for $\gamma$ is too small, we encounter overflow.
    The two's complement mantissa then wraps around, producing large errors in the result.
    This can be circumvented by saturation of the result to the range $[-2^{\out{w}-1}, 2^{\out{w}-1} - 1]$ of the
    two's complement integer.
    Then again, small deviations of $\gamma$ from the actual infinity norm of the result are only leading to small
    errors in the computed vector.

    These considerations are manifested in a modified version of \cref{alg:qcomp} (\texttt{qcomp}) listed in
    \cref{alg:qcomp-saturate} (\texttt{nnqcomp}).
    It skips the normalization step, and already returns after the first loop.
    The atomic update in \cref{alg:qcomp}, line~\ref{alg:qcomp:line:atomic}, is removed, the computed vector entries are
    saturated before being truncated to the target precision, and no temporary vector is required, saving the additional memory.
    \noindent
\begin{algorithm}[t]
    \begin{center}
        \begin{minipage}{\linewidth}
            \footnotesize
            \caption{$\texttt{Non-normalized quantized BFP vector computation: } \texttt{nnqcomp}$}
            \label{alg:qcomp-saturate}
            \begin{algorithmic}[1]

                % Input arguments
                \Require $
                \texttt{setup-func},\
                \texttt{comp-func},\
                \texttt{input},\
                \out{w}: \unsignedpositive,\
                \gamma: \bfp^1,\
                n: \unsignedpositive$

                % \Ensure $0 < \out{w} \leq \tmp{w} \leq \bfpq{\exact{z}}, \ \gamma > 0,\ \gamma \text{ normalized}$
                \Ensure $0 < \out{w}, \ \gamma > 0,\ \gamma \text{ normalized}$

                \Variables $
                \texttt{setup},\
                \exact{z}: \bfp^1\ (\text{exact result}), \
                \out{z}:\bfp^n\ (\text{result}),\
                % (\texttt{output-setup}),\
                \out{\mu},\
                % \exact{\mu},\
                \lambda_{\textrm{out}} : \signed$

                \LineComment{Allocation and setup.}

                % \State $\exact{z},\ (\texttt{output-setup}) \gets \texttt{setup-func}(\texttt{input})$
                \State $\left(\exact{z},\ \texttt{setup}\right)  \gets \texttt{setup-func}(\texttt{input})$
                \Comment{Setup necessary variables for exact computations.}

                \State $\out{\mu} \gets \texttt{msb}(\bfpm{\gamma}) + \bfpe{\gamma} - \bfpe{\exact{z}}$
                \Comment{MSB of the estimated mantissa window.}

                \State $\out{\lambda} \gets \out{\mu} - \out{w}$
                \Comment{LSB of the estimated mantissa window.}

                \State $\bfpq{(\out{z})} \gets \out{w}$
                \Comment{Width of result BFP vector.}

                \State $\bfpe{(\out{z})} \gets \bfpe{\exact{z}} + \out{\lambda}$
                \Comment{Block-exponent of truncated result $\out{z}$.}

                \LineComment{Computation of the saturated and truncated result BFP vector.}

                %\State $\exact{\mu} \gets 1$
                %\Comment{Keep track of the largest MSB idx. of the exact result.}

                \ParFor{$i = 1,\, \dots,\, n$}

                \State $\exact{z} \gets \texttt{comp-func}(\exact{z},\ \texttt{input},\ \texttt{setup},\ i)$
                %\label{alg:qcomp:line:firstloopexact}

                \Comment{Compute exact quantity in parallel.}

                \State $\bfpm{\exact{z}} \gets \bfpm{\exact{z}} \RightShift \out{\lambda}$
                \Comment{Shift and truncate exact result.}

                \State $\bfpm{\exact{z}} \gets \texttt{clamp}\left( \bfpm{\exact{z}}, -2^{\out{w} - 1}, 2^{\out{w}-1} - 1 \right)$

                \Comment{Saturate result.}

                \State $(\bfpm{(\out{z})})_i \gets \text{decr} \left( \bfpq{\exact{z}} - \bfpq{(\out{z})},\ \bfpm{\exact{z}} \right)$
                %\label{alg:qcomp:line:firstlooptrunc}

                \Comment{Cast result into output BFP vector.}

                \EndParFor

                \State\Return $\out{z}$
                \label{alg:qcomp:line:early-return}
                \Comment{The result is saturated and truncated.}

            \end{algorithmic}
        \end{minipage}
    \end{center}
\end{algorithm}

    \begin{figure}
        \centering
        \resizebox*{\textwidth}{!}{%
            \input{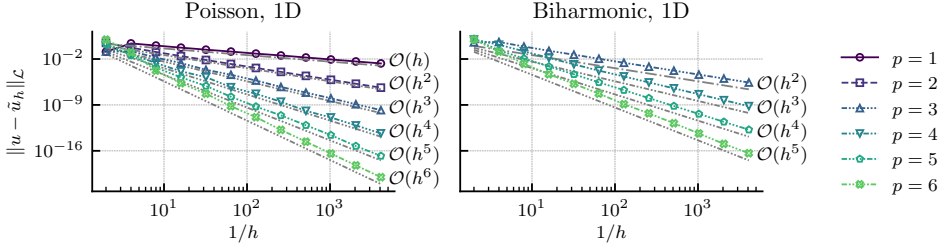}
        }
        \vspace{-24pt}
        \caption{Progressive precision \gls*{bfp}-\fmg~using the same setup as in \cref{fig:bfp-fmg} but without
        \gls*{bfp} vector normalization, i.e., \cref{alg:qcomp} (\texttt{qcomp}) is replaced by \cref{alg:qcomp-saturate} (\texttt{nnqcomp})
        for all vector operations on all refinement levels.
        Still, the solver produces discretization-error-accurate approximations in all tested cases.
        There is no significant error difference compared to the results computed with \cref{alg:qcomp}.}
        \label{fig:bfp-fmg-saturated}
    \end{figure}

    \Cref{fig:bfp-fmg-saturated} shows the total energy error over refinement using progressive precision \gls*{bfp}-\fmg~
    with \cref{alg:qcomp-saturate} instead of \cref{alg:qcomp}.
    The heuristic estimates for $\gamma$ are sufficient to ensure discretization-error-accuracy for both model problems,
    even without \gls*{bfp}-vector normalization.
    Direct comparison with the results using the safe \cref{alg:qcomp} shows no significant differences in error.
    We have briefly tested this also for higher order discretizations ($p = 10$), observing the same results.

    \begin{remark}[Skipping normalization in \ir-\V]\label{rem:skipping-normalization-ir-v}
        The situation is slightly different for \ir-\V~applied directly (i.e., without \fmg) to a zero initial guess.
        During the first few (1-2) iterations, we observe that the infinity-norm of the \ir-residual (left-hand side in
        \cref{alg:ir}, line~\ref{alg:ir-residual}) may in fact increase before it decreases.
        This increase may be so rapid that the estimates of $\gamma$ from \cref{tab:estimates} yield large overflows,
        leading to a diverging iteration. (This is especially pronounced for the lowest-order approximations, i.e.,
        Poisson $p = 1$, and biharmonic $p = 3$; higher-order approximations are less problematic.)
        Using \texttt{qcomp} (\cref{alg:qcomp}) instead of \texttt{nnqcomp} (\cref{alg:qcomp-saturate}) only for the
        \ir-residual computation (\cref{alg:ir}, line~\ref{alg:ir-residual}) during the first few iterations of
        \ir-\V~fixes this issue, and we observe convergence for all test cases, if normalization is skipped for
        all remaining \gls*{bfp}-operations.
        Due to the initial approximation from the coarse grids, this is not necessary inside of \fmg.
    \end{remark}

    Skipping \gls*{bfp}-vector normalization promises significant performance advantages in practice, without compromise
    regarding solution accuracy during \fmg.
    An iterative solver using \cref{alg:qcomp-saturate} can additionally be equipped with a safety mechanism, that
    observes the residual convergence.
    Similar to the approach in \cref{rem:skipping-normalization-ir-v}, \cref{alg:qcomp} can then be invoked
    dynamically, for certain operations, and for a subset of iterations, if the convergence stalls.
    However, as shown in \cref{fig:bfp-fmg-saturated}, we do not observe that this is necessary for
    \gls*{bfp}-\fmg~applied to our model problems.

    \section{Conclusion}\label{sec:conclusion}

    This paper has demonstrated that the solution of elliptic \glspl*{pde} in pure \emph{integer arithmetic} can be done in
    practice.
    The results lay the groundwork for energy efficient implementations on specialized hardware.
    Additionally, the asymptotic precision bounds
    from~\cite{McCormick:2021:AlgebraicErrorAnalysis,Tamstorf:2021:DiscretizationErrorAccurateMixedPrecisionMultigrid}
    have been applied successfully to obtain discretization-error-accuracy using a progressive- and mixed-precision
    multigrid solver in \gls*{bfp} format.
    To achieve this, we have proposed an efficient \gls*{bfp}-algorithm that ensures exact computation of common BLAS-like vector-vector
    and matrix-vector operations up to a specified precision.
    Compared to the results in~\cite[Figures 2, 9]{Tamstorf:2021:DiscretizationErrorAccurateMixedPrecisionMultigrid}, the number of
    bits required to attain a certain level of error accuracy is similar to that of standard floating point arithmetic.
    \gls*{bfp} arithmetic is particularly efficient when applied to \fmg, since our results suggest that normalization of the \gls*{bfp} vectors is not necessary
    with proper estimates of upper bounds of the infinity norm of the intermediate results.
    Hopefully, all of this will stimulate future research in this area to establish
    a rigorous theoretical framework for iterative linear solvers in \gls*{bfp} arithmetic, and to develop
    accurate performance models with an eye towards deployment of these ideas for real applications. 

    \bibliographystyle{siamplain}
    \bibliography{references}

\begin{thebibliography}{10}

\bibitem{BasicLinearAlgebraSubprogramsTechnicalForum:2002:BasicLinearAlgebra}
{\sc {Basic Linear Algebra Subprograms Technical Forum}}, {\em Basic {{Linear
  Algebra Subprograms Technical Forum Standard}}}, International Journal of
  High Performance Applications and Supercomputing, 16 (2002), pp.~1--111,
  \url{https://journals.sagepub.com/toc/hpcc/16/1}.

\bibitem{Basumallik2022}
{\sc A.~Basumallik, D.~Bunandar, N.~Dronen, N.~Harris, L.~Levkova, C.~McCarter,
  L.~Nair, D.~Walter, and D.~Widemann}, {\em {Adaptive Block Floating-Point for
  Analog Deep Learning Hardware}}, 2022,
  \url{https://doi.org/10.48550/arXiv.2205.06287}.
\newblock Under submission at IEEE Transactions on Neural Networks and Learning
  Systems (TNNLS).

\bibitem{Boldo2020}
{\sc S.~Boldo, D.~Gallois-Wong, and T.~Hilaire}, {\em {A Correctly-Rounded
  Fixed-Point-Arithmetic Dot-Product Algorithm}}, in Proceedings of the 27th
  IEEE Symposium on Computer Arithmetic, ARITH-2020, IEEE Computer Society,
  June 2020, pp.~9--16, \url{https://doi.org/10.1109/ARITH48897.2020.00011}.

\bibitem{Dai2021}
{\sc S.~Dai, R.~Venkatesan, H.~Ren, B.~Zimmer, W.~J. Dally, and B.~Khailany},
  {\em {VS-Quant: Per-vector Scaled Quantization for Accurate Low-Precision
  Neural Network Inference}}, in Proceedings of Machine Learning and Systems,
  A.~Smola, A.~Dimakis, and I.~Stoica, eds., vol.~3 of MLSys, 2021,
  pp.~873--884,
  \url{https://proceedings.mlsys.org/paper/2021/file/f0935e4cd5920aa6c7c996a5ee53a70f-Paper.pdf}.

\bibitem{Drumond2018}
{\sc M.~Drumond, T.~Lin, M.~Jaggi, and B.~Falsafi}, {\em {Training DNNs with
  Hybrid Block Floating Point}}, in Proceedings of NeurIPS'18, S.~Bengio,
  H.~Wallach, H.~Larochelle, K.~Grauman, N.~Cesa-Bianchi, and R.~Garnett, eds.,
  vol.~31 of Advances in Neural Information Processing Systems, Red Hook, NY,
  USA, 2018, Curran Associates Inc., p.~453–463,
  \url{https://proceedings.neurips.cc/paper/2018/file/6a9aeddfc689c1d0e3b9ccc3ab651bc5-Paper.pdf}.

\bibitem{Gustafson2017}
{\sc Gustafson and Yonemoto}, {\em {Beating Floating Point at Its Own Game:
  Posit Arithmetic}}, Supercomputing Frontiers and Innovations: an
  International Journal, 4 (2017), p.~71–86,
  \url{https://doi.org/10.14529/jsfi170206}.

\bibitem{Gutknecht2002}
{\sc M.~H. Gutknecht and S.~R{\"o}llin}, {\em {The Chebyshev iteration
  revisited}}, Parallel Computing, 28 (2002), pp.~263--283,
  \url{https://doi.org/10.1016/S0167-8191(01)00139-9}.

\bibitem{Horowitz2014}
{\sc M.~Horowitz}, {\em {Computing's Energy Problem (and what we can do about
  it)}}, in 2014 IEEE International Solid-State Circuits Conference Digest of
  Technical Papers (ISSCC), IEEE, Feb. 2014, pp.~10--14,
  \url{https://doi.org/10.1109/ISSCC.2014.6757323}.

\bibitem{Jouppi2021}
{\sc N.~P. Jouppi, D.~Hyun~Yoon, M.~Ashcraft, M.~Gottscho, T.~B. Jablin,
  G.~Kurian, J.~Laudon, S.~Li, P.~Ma, X.~Ma, T.~Norrie, N.~Patil, S.~Prasad,
  C.~Young, Z.~Zhou, and D.~Patterson}, {\em {Ten Lessons From Three
  Generations Shaped Google’s TPUv4i : Industrial Product}}, in 2021 ACM/IEEE
  48th Annual International Symposium on Computer Architecture (ISCA), 2021,
  pp.~1--14, \url{https://doi.org/10.1109/ISCA52012.2021.00010}.

\bibitem{Koster2017}
{\sc U.~K\"{o}ster, T.~J. Webb, X.~Wang, M.~Nassar, A.~K. Bansal, W.~H.
  Constable, O.~H. Elibol, S.~Gray, S.~Hall, L.~Hornof, A.~Khosrowshahi,
  C.~Kloss, R.~J. Pai, and N.~Rao}, {\em {Flexpoint: An Adaptive Numerical
  Format for Efficient Training of Deep Neural Networks}}, in Proceedings of
  NIPS'17, I.~Guyon, U.~V. Luxburg, S.~Bengio, H.~Wallach, R.~Fergus,
  S.~Vishwanathan, and R.~Garnett, eds., vol.~30 of Advances in Neural
  Information Processing Systems, Red Hook, NY, USA, 2017, Curran Associates
  Inc., p.~1742–1752,
  \url{https://proceedings.neurips.cc/paper/2017/file/a0160709701140704575d499c997b6ca-Paper.pdf}.

\bibitem{Kulisch2011b}
{\sc U.~Kulisch}, {\em {Very fast and exact accumulation of products}},
  Computing, 91 (2011), pp.~397--405,
  \url{https://doi.org/10.1007/s00607-010-0131-y}.

\bibitem{Lian2019}
{\sc X.~Lian, Z.~Liu, Z.~Song, J.~Dai, W.~Zhou, and X.~Ji}, {\em
  {High-Performance FPGA-Based CNN Accelerator With Block-Floating-Point
  Arithmetic}}, IEEE Transactions on Very Large Scale Integration (VLSI)
  Systems, 27 (2019), pp.~1874--1885,
  \url{https://doi.org/10.1109/TVLSI.2019.2913958}.

\bibitem{McCormick:2021:AlgebraicErrorAnalysis}
{\sc S.~F. McCormick, J.~Benzaken, and R.~Tamstorf}, {\em {Algebraic Error
  Analysis for Mixed-Precision Multigrid Solvers}}, SIAM Journal on Scientific
  Computing, 43 (2021), pp.~S392--S419,
  \url{https://doi.org/10.1137/20M1348571}.

\bibitem{Noh2022a}
{\sc S.-H. Noh, J.~Koo, S.~Lee, J.~Park, and J.~Kung}, {\em {FlexBlock: A
  Flexible DNN Training Accelerator with Multi-Mode Block Floating Point
  Support}}, 2022, \url{https://doi.org/10.48550/arXiv.2203.06673}.
\newblock Under revision at IEEE Transactions on Computers.

\bibitem{Noh2022b}
{\sc S.-H. Noh, J.~Park, D.~Park, J.~Koo, J.~Choi, and J.~Kung}, {\em
  {LightNorm: Area and Energy-Efficient Batch Normalization Hardware for
  On-Device DNN Training}}, 2022,
  \url{https://doi.org/10.48550/arXiv.2211.02686}.

\bibitem{Behrooz:2010}
{\sc B.~Parhami}, {\em {Computer Arithmetic: Algorithms and Hardware Designs}},
  Oxford University Press, New York, 2nd~ed., 2010.

\bibitem{QianZhang:2022:FASTDNNTraining}
{\sc S.~Qian~Zhang, B.~McDanel, and H.~T. Kung}, {\em {FAST: DNN Training Under
  Variable Precision Block Floating Point with Stochastic Rounding}}, in 2022
  IEEE International Symposium on High-Performance Computer Architecture
  (HPCA), IEEE, Apr. 2022, pp.~846--860,
  \url{https://doi.org/10.1109/HPCA53966.2022.00067}.

\bibitem{Rouhani2020}
{\sc B.~D. Rouhani, D.~Lo, R.~Zhao, M.~Liu, J.~Fowers, K.~Ovtcharov,
  A.~Vinogradsky, S.~Massengill, L.~Yang, R.~Bittner, A.~Forin, H.~Zhu, T.~Na,
  P.~Patel, S.~Che, L.~Chand~Koppaka, X.~Song, S.~Som, K.~Das, S.~Tiwary,
  S.~Reinhardt, S.~Lanka, E.~Chung, and D.~Burger}, {\em {Pushing the Limits of
  Narrow Precision Inferencing at Cloud Scale with Microsoft Floating Point}},
  in Proceedings of NeurIPS 2020, H.~Larochelle, M.~Ranzato, R.~Hadsell,
  M.~Balcan, and H.~Lin, eds., vol.~33 of Advances in Neural Information
  Processing Systems, Curran Associates, Inc., November 2020, pp.~10271--10281,
  \url{https://proceedings.neurips.cc/paper/2020/hash/747e32ab0fea7fbd2ad9ec03daa3f840-Abstract.html}.

\bibitem{Tamstorf:2021:DiscretizationErrorAccurateMixedPrecisionMultigrid}
{\sc R.~Tamstorf, J.~Benzaken, and S.~F. McCormick}, {\em
  {Discretization-Error-Accurate Mixed-Precision Multigrid Solvers}}, SIAM
  Journal on Scientific Computing, 43 (2021), pp.~S420--S447,
  \url{https://doi.org/10.1137/20M1349230}.

\bibitem{Wilkinson1963}
{\sc J.~H. Wilkinson}, {\em {Rounding Errors in Algebraic Processes}},
  Prentice-Hall series in Automatic Computation, Prentice-Hall, Englewood
  Cliffs, N. J., 1963.

\end{thebibliography}

    \appendix
 
     \section{Discrete Harmonic}
    \label{disc-harm}
   Using the terminology in Section~\ref{sec:bfp-quant-error}, suppose that $v$ has minimum energy, $\|v\|_A$, subject to the constraint $\|v\|_\infty = 1$. To estimate $\|v\|_A$, assume without loss of generality that $v$ is $1$ at grid point $p$: $v_p = 1$. Note that $v$ is a discrete harmonic in the sense that $(Av)_q = 0$ for all grid points $q \ne p$. To see this, letting $s$ be a scalar and $d \equiv (Av)_q$ for any $q \ne p$, we would then have that
\[
\langle A(v - sd), v - sd \rangle = \langle Av, v \rangle - 2s \langle Av, d \rangle + s^2 \langle Ad, d \rangle
= \langle Av, v \rangle - 2s \|d\|^2 + s^2 \langle Ad, d \rangle .
\]
If $d$ were not $0$, then choosing $s > 0$ small enough (e.g., $s < \frac{\langle Ad, d \rangle}{\|d\|^2}$) would mean that $\langle A(v - sd), v - sd \rangle < \langle Av, v \rangle$,
which contradicts optimality of $v$. Hence, $v = \gamma A^{-1}e_p$, where $\gamma = \frac{1}{\langle e_p, A^{-1} e_p \rangle}$ and $e_p$ is the vector that is $1$ at $p$ and $0$ elsewhere. (Note that $v$ as defined here satisfies $v_p = \langle e_p, v \rangle = \gamma \langle e_p, A^{-1}e_p \rangle = 1$ and $v_q = \langle e_q, Av \rangle = \gamma \langle e_q, e_p \rangle = 0$ for $q \ne p$.) Thus,
\[
\|v\|_A = \langle Av, v \rangle^\frac{1}{2} = \gamma \langle e_p, A^{-1}e_p \rangle^\frac{1}{2} = \gamma^\frac{1}{2} = \langle e_p, A^{-1}e_p\rangle^{-\frac{1}{2}} .
\]
We therefore have that the minimum value of $\|v\|_A$ is the inverse square root of the maximum diagonal entry of $A^{-1}$.
 
    \section{BFP BLAS algorithms}
    \label{sec:appendix-bfp-blas-algos}

    \begin{algorithm}[H]
    \begin{center}
        \begin{minipage}{\linewidth}
            \footnotesize
            \caption{$\texttt{Exact axpby setup: }
            \texttt{eaxpby-setup}$}
            \label{alg:eaxpby-setup}
            \begin{algorithmic}[1]

                % Input arguments
                \Require $
                x : \bfp^{n}, \
                y : \bfp^{n}, \
                \alpha : \bfp^{1}, \, \
                \beta : \bfp^{1}
                $

                % Local variables
                \Variables $
                \exact{a} : \bfp^{1}, \
                \exact{b} : \bfp^{1}, \
                \exact{z} : \bfp^{1}, \
                d : \signed
                $

                \LineComment{Setup for the two products $\exact{a} = \alpha x$ and $\exact{b} = \beta y$.}

                \State $\bfpq{\exact{a}} \gets \bfpq{\alpha} + \bfpq{x}$
                \Comment{Widths of the products.}

                \State $\bfpq{\exact{b}} \gets \bfpq{\beta} + \bfpq{y}$

                \State $\bfpe{\exact{a}} \gets \bfpe{\alpha} + \bfpe{x}$
                \Comment{Block-exponents of the products.}

                \State $\bfpe{\exact{b}} \gets \bfpe{\beta} + \bfpe{y}$

                \LineComment{Setup for the sum $z = \exact{a} + \exact{b}$.}

                \State $d \gets \bfpe{\exact{a}} - \bfpe{\exact{b}}$
                \Comment{Difference of block exponents of $\alpha x$ and $\beta y$.}

                \If{$(d < 0)$}
                    \Comment{Aligning block-exponents of $\alpha x$ and $\beta y$.}

                    \State $\bfpq{\exact{z}} \gets \max( \bfpq{\exact{a}}, \bfpq{\exact{b}} + |d| ) + 1$
                \Else
                    \State $\bfpq{\exact{z}} \gets \max( \bfpq{\exact{b}}, \bfpq{\exact{a}} + |d| ) + 1$
                \EndIf

                \State $\bfpe{\exact{z}} \gets \min( \bfpe{\exact{a}}, \bfpe{\exact{b}} )$
                \Comment{Block-exponent of $\exact{z} = \alpha x + \beta y$.}

                \State\Return $\left( \exact{z},\ (\exact{a}, \exact{b}, d) \right)$

            \end{algorithmic}
        \end{minipage}
    \end{center}
\end{algorithm}

    \begin{algorithm}[H]
    \begin{center}
        \begin{minipage}{\linewidth}
            \footnotesize
            \caption{$\texttt{Exact axpby row: }
            \texttt{eaxpby-row}$}
            \label{alg:eaxpby-row}
            \begin{algorithmic}[1]

                % Input arguments
                \Require $
                \exact{z} : \bfp^{1}, \
                (
                x : \bfp^{n}, \
                y : \bfp^{n}, \
                \alpha : \bfp^{1}, \
                \beta : \bfp^{1}
                ),\
                (
                \exact{a} : \bfp^{1}, \
                \exact{b} : \bfp^{1}, \
                d : \signed
                ), \
                i: \unsignedpositive$

                % Local variables
                % \Variables $$

                % ----- multiplication ----- %

                \LineComment{Exact scalar multiplication.}

                \State $\bfpm{\exact{a}} \gets \bfpm{\alpha} \cdot (\bfpm{x})_i$
                \State $\bfpm{\exact{b}} \gets \bfpm{\beta}  \cdot (\bfpm{y})_i$

                % ----- ADD ----- %

                \LineComment{Exact addition. Requires alignment of block-exponents of $\exact{a}$ and $\exact{b}$.}

                \If{($d < 0$)}
                    \Comment{$\bfpe{\exact{a}} < \bfpe{\exact{b}}$, resulting block-exponent is $\bfpe{\exact{a}}$.}

                    \State $\bfpm{\exact{b}} \gets \text{incr}\left(|d|, \bfpm{\exact{b}} \right)$
                    \Comment{Ensure enough space to left-shift.}

                    \State $\bfpm{\exact{b}} \gets \bfpm{\exact{b}} \LeftShift |d|$
                    \Comment{Shifting mantissa to align exponents of $\exact{a}$ and $\exact{b}$.}

                \Else
                    \Comment{$\bfpe{\exact{a}} \geq \bfpe{\exact{b}}$, resulting block-exponent is $\bfpe{\exact{b}}$.}

                    \State $\bfpm{\exact{a}} \gets \text{incr}\left(|d|, \bfpm{\exact{a}} \right)$
                    \Comment{Ensure enough space to left-shift.}

                    \State $\bfpm{\exact{a}} \gets \bfpm{\exact{a}} \LeftShift |d|$
                    \Comment{Shifting mantissa to align exponents of $\exact{a}$ and $\exact{b}$.}

                \EndIf

                \State $\bfpm{\exact{z}} \gets \bfpm{\exact{a}} + \bfpm{\exact{b}}$
                \Comment{Computes exact sum $\bfpm{\exact{z}} = \bfpm{\exact{a}} + \bfpm{\exact{b}}$.}

                % ----- Return ----- %

                \State\Return $\exact{z}$
                \Comment{Return row $i$ of the exact $\alpha x + \beta y$.}

            \end{algorithmic}
        \end{minipage}
    \end{center}
\end{algorithm}

    \begin{algorithm}[H]
    \begin{center}
        \begin{minipage}{\linewidth}
            \footnotesize
            \caption{$\texttt{Exact SpMV setup: }
            \texttt{espmv-setup}$}
            \label{alg:ematvec-setup}
            \begin{algorithmic}[1]

                % Input arguments
                \Require $
                A : \bfp^{n \times m}, \
                x: \bfp^{m}, \
                m_A: \unsignedpositive \ \text{(max. num. non-zeros per row of $A$)}
                $

                % Local variables
                \Variables $
                \exact{z} : \bfp^1
                $

                \LineComment{Setup for the exact product $Ax$.}

                \State $\bfpq{\exact{z}} \gets \bfpq{A} + \bfpq{x} + \lceil \log_2(m_A) \rceil$
                \Comment{Width of mantissa of exact $Ax$.}

                \State $\bfpe{\exact{z}} \gets \bfpe{A} + \bfpe{x}$
                \Comment{Block-exponent of exact $Ax$.}

                \State\Return $\left( \exact{z},\ () \right)$

            \end{algorithmic}
        \end{minipage}
    \end{center}
\end{algorithm}

    \begin{algorithm}[H]
    \begin{center}
        \begin{minipage}{\linewidth}

            \footnotesize
            \caption{$\texttt{Exact SpMV row: }
            \texttt{espmv-row}$}
            \label{alg:ematvec}
            \begin{algorithmic}[1]

                % Input arguments
                \Require $
                \exact{z}: \bfp^1, \
                (
                A : \bfp^{n \times m},\
                x : \bfp^{m},\
                m_A : \unsignedpositive
                ),\
                (),\
                i : \unsignedpositive
                $

                % Local variables
                \Variables $
                \exact{t} : \twoscomp{\bfpq{A} + \bfpq{x}}
                % \bfpm{z} : \textcolor{dtype}{\twoscomp{\bfpq{z}}},\
                % u : \textcolor{dtype}{\twoscomp{\bfpq{z}}}, \
                % g : \textcolor{dtype}{\twoscomp{\bfpq{b}}}, \
                %\bfpm{\exact{r}} : \textcolor{dtype}{\twoscomp{\bfpq{\exact{r}}}},
                %\bfpm{\tilde{r}} : \textcolor{dtype}{\twoscomp{\bfpq{\tilde{r}}}}
                $

                \LineComment{Exact dot product of one row of $A$ and $x$.}

                \State $\bfpm{\exact{z}} \gets 0$

                \For{$j \in \{ [1, n] \cap \unsignedpositive : (\bfpm{\mathbf{A}})_{ij} \neq 0 \}$}
                    \State $\exact{t} \gets (\bfpm{\mathbf{A}})_{ij} \cdot (\bfpm{\mathbf{x}})_j$
                    \Comment{Exact integer multiplications.}

                    \State $\bfpm{\exact{z}} \gets \bfpm{\exact{z}} + \exact{t}$
                    \Comment{Exact integer accumulation.}
                \EndFor

                % ----- Return ----- %

                \State\Return $\exact{z}$
                \Comment{Return row $i$ of the exact $Ax$.}

            \end{algorithmic}
        \end{minipage}
    \end{center}
\end{algorithm}

    \begin{algorithm}[H]
    \begin{center}
        \begin{minipage}{\linewidth}
            \footnotesize
            \caption{$\texttt{Exact gemv setup: }
            \texttt{egemv-setup}$}
            \label{alg:egemv-setup}
            \begin{algorithmic}[1]

                % Input arguments
                \Require $
                A : \bfp^{n \times m},\
                x : \bfp^{m},\
                y : \bfp^{n},\
                \alpha : \bfp^{1},\
                \beta : \bfp^{1},\
                m_A : \unsignedpositive
                $

                % Local variables
                \Variables $
                \exact{g} : \bfp^{1},\
                \exact{z} : \bfp^{1},\
                $

                \LineComment{Reusing \texttt{espmv} and \texttt{eaxpby}.}

                \State $\left(\exact{g},\ \texttt{<empty>}\right) \gets \texttt{espmv-setup}(A,\ x,\ m_A)$

                \Comment{Setup for $\exact{g} = Ax$.}

                \State $\left(\exact{z},\ (\exact{a}, \exact{b}, d)\right) \gets \texttt{eaxpby-setup}(\exact{g},\ y,\ \alpha,\ \beta)$

                \Comment{Setup for $\alpha \exact{g} + \beta y$.}

                \State\Return $\left( \exact{z},\ (\exact{g}, \exact{a}, \exact{b}, d) \right)$

            \end{algorithmic}
        \end{minipage}
    \end{center}
\end{algorithm}

    \begin{algorithm}[H]
    \begin{center}
        \begin{minipage}{\linewidth}
            \footnotesize
            \caption{$\texttt{Exact gemv row: }
            \texttt{egemv-row}$}
            \label{alg:egemv-row}
            \begin{algorithmic}[1]

                % Input arguments
                \Require $
                \exact{z} : \bfp^{1}, \
                (
                A : \bfp^{n \times m},\
                x : \bfp^{m},\
                y : \bfp^{n},\
                \alpha : \bfp^{1},\
                \beta : \bfp^{1},\
                m_A : \unsignedpositive
                ), \newline
                (
                \exact{g} : \bfp^{1},\
                \exact{a} : \bfp^{1},\
                \exact{b} : \bfp^{1},\
                d : \unsigned
                ),\
                i : \unsignedpositive
                $

                % Local variables
                % \Variables $$

                \LineComment{Reusing \texttt{espmv} and \texttt{eaxpby}.}

                \State $\exact{g} \gets \texttt{espmv-row}(\exact{g},\ (A, x, m_A),\ \texttt{<empty>},\ i)$

                \Comment{Computing one row of $Ax$.}

                \State $\exact{z} \gets \texttt{eaxpby-row}(\exact{z}, (\exact{g}, y_i, \alpha, \beta), (\exact{a}, \exact{b}, d), 1)$

                \Comment{Computing one row of $\alpha Ax + \beta y$.}

                % ----- Return ----- %

                \State\Return $\exact{z}$
                \Comment{Return row $i$ of the exact $\alpha Ax + \beta y$.}

            \end{algorithmic}
        \end{minipage}
    \end{center}
\end{algorithm}

    \section{BFP multigrid algorithms}
    \label{sec:appendix-bfp-mg-algos}

    \begin{algorithm}[H]
    \begin{center}
        \begin{minipage}{1\linewidth}
            \footnotesize
            \caption{$\texttt{Coefficients for two Chebyshev iterations}$}
            \label{alg:cheby-nodes}
            \begin{algorithmic}[1]
                \Require $\rho$ (upper bound for max. eigenvalue of generalized problem $Ax = \lambda Dx$), \newline
                $0 < \eta < 1$ (part of spectrum to target)
                \State $\alpha \leftarrow \tfrac{1}{2}(1 + \eta)\rho$
                \State $c \leftarrow \tfrac{1}{2}(1-\eta)\rho$
                \State $\beta \leftarrow \alpha -\frac{c^2}{2\alpha}$
                \State $c_1 \leftarrow 2/\beta$
                \State $c_2 \leftarrow -1/(\alpha\beta)$
                \State\Return $(c_1,c_2)$\COMMENT{Return coefficients}
            \end{algorithmic}
        \end{minipage}
    \end{center}
\end{algorithm}

    \begin{algorithm}[H]
    \begin{center}
        \begin{minipage}{1\linewidth}
            \footnotesize
            \caption{$\texttt{Setup ({\s}) to ensure that }D=I$}
            \label{alg:solver-setup}
            \begin{algorithmic}[1]
                \Require $A_i$, $b_i$, $1 \le i \le \ell$, $P_j, 2 \le j \le \ell$ (prolongation operators), $\ell \ge 1$ (number of levels), \newline
                $\ell_{\mathrm{est}}$ (level on which the Chebyshev nodes are estimated),\newline
                $\mu^*$ (part of the spectrum to target with the Chebyshev smoother)
                \State {$\rho \leftarrow \texttt{MaxGenEigenvalueUpperBound}(A_{\ell_{\mathrm{est}}}, D_{\ell_{\mathrm{est}}})$}
                \State {$(c_{1}, c_{2}) \leftarrow \texttt{ChebyshevNodes}(\rho, \mu^*)$}
                \State $i \leftarrow \ell$\COMMENT{Initialize {\s } }
                \State $D_f \leftarrow \diag(A_i)$
                \While {$i > 0$}
                    \State $A_i \leftarrow D_f^{-1}A_i$ \COMMENT{Premultiply $A$ by $D$ so new $D=I$}
                    \State $b_i \leftarrow D_f^{-1}b_i$ \COMMENT{Preserve solution of fine-level equation}
                    \If{$i>1$}
                        \State $D_c \leftarrow \diag(A_{i-1})$
                        \State $R_i \leftarrow D_c^{-1}P_i^T D_f$\COMMENT{Preserve solution of correction equation}
                        \State $D_f \leftarrow D_c$
                    \EndIf
                    \State $i \leftarrow i - 1$\COMMENT{Decrement {\s } cycle counter}
                \EndWhile
                % \State\Return $A_i$, $b_i$, $1 \le i \le \ell$, $R_j$, $2 \le j \le \ell$, $c_{1}$, $c_{2}$
                \State\Return $\left( A_1, \dots, A_\ell, b_1, \dots, b_\ell, R_2, \dots, R_\ell, c_1, c_2 \right)$
            \end{algorithmic}
        \end{minipage}
    \end{center}
\end{algorithm}

    \vspace{-14pt}
    \begin{algorithm}[H]
    \begin{center}
        \begin{minipage}{1\linewidth}
            \footnotesize
            \caption{$\texttt{Iterative Refinement (\ir) with }D=I$}
            \label{alg:ir}
            \begin{algorithmic}[1]
                \Require $A$, $b$, $x$ (initial guess), tol $> 0$ (convergence tolerance), $m_A$, \texttt{InnerSolver}

                \State $\left(A, b\right) \gets \left(\texttt{quant}(A), \texttt{quant}(b)\right)$
                \COMMENT{Quantize to $\weweq$ bits}

                %\State $b \gets b$
                %\COMMENT{Quantize {$b$} to $\weweq$ bits}

                \State {{$r \leftarrow$}$\texttt{qgemv}(A, x, b, 1, -1, m_A, \wewed)$}
                \COMMENT{{Compute {\ir } Residual $r \leftarrow$}$Ax - b$}\label{alg:ir-residual}

                \If{$\|r\| <$ tol}
                    \State\Return $x$
                    \COMMENT{Return Solution of $Ax = b$}
                \EndIf
                \State $y \leftarrow \texttt{InnerSolver}(A, r)$
                \COMMENT{Approximate Solution $y$ of $Ay = r$}\label{ir-inner}

                \State $x \leftarrow \texttt{qsub}(x, y, \wewe)$
                \COMMENT{Update Approximation $x \leftarrow x - y$}\label{alg:ir-correction}
                \State\Goto{alg:ir-residual}
            \end{algorithmic}
        \end{minipage}
    \end{center}
\end{algorithm}

    \vspace{-14pt}
    \begin{algorithm}[H]
    \begin{center}
        \begin{minipage}{1\linewidth}
            \footnotesize
            \caption{$\texttt{V(1,0)-Cycle (\V) Correction Scheme with }D=I$}
            \label{alg:v-cycle}
            \begin{algorithmic}[1]
                \Require $A, r, P, R, \ell \ge 1$ (number of {\V } levels); $c_1, c_2$ (Chebyshev coefficients), $m_A, m_P, m_R$

                \State $\left(A, P, R\right) \leftarrow \left(\texttt{quant}(A), \texttt{quant}(P), \texttt{quant}(R)\right)$
                \COMMENT{Quantize to $\wewed_\ell$ bits}

                \State {$y \leftarrow \texttt{qgemv}(A, r, r, c_2, c_1, m_A, \wewed_\ell)$}\label{alg:v-cycle-relaxation}
                \COMMENT{{Relax on Current Approximation ($y=0$)}}

                \If {$\ell > 1$}
                \COMMENT{Check for Coarser Level}

                    \State {$r_{\textrm{v}} \leftarrow \texttt{qgemv}(A, y, r, 1, -1, m_A, \wewed_\ell)$}\label{alg:v-cycle-residual}
                    \COMMENT{{Evaluate {\V } Residual $r_{\textrm{v}} \leftarrow A y - r$}}

                    \State {$r_{\ell - 1} \leftarrow \texttt{qspmv}(R, r_{\textrm{v}}, m_R, \wewed_\ell)$}\label{alg:v-cycle-restriction}
                    \COMMENT{{Restrict {\V } Residual $r_{\ell - 1} \leftarrow Rr_{\textrm{v}}$}}

                    \State {$d_{\ell - 1} \leftarrow $\V$(A_{\ell - 1}, r_{\ell - 1}, P_{\ell - 1}, R_{\ell - 1}, \ell - 1, m_A, m_P, m_R)$}\\
                    \COMMENT{{Compute Correction from Coarser Levels}}\label{v-cgc}

                    \State {$y \leftarrow \texttt{qgemv}(P, d_{\ell - 1}, y, -1, 1, m_P, \wewed_\ell)$}\label{alg:v-cycle-correction}
                    \COMMENT{{Interpolate \& Update $y \leftarrow y - Pd_{\ell - 1}$}}

                \EndIf
                \State\Return ${y}$\COMMENT{{Return Approximate Solution of $A y = r$}}
            \end{algorithmic}
        \end{minipage}
    \end{center}
\end{algorithm}

    \vspace{-14pt}
    \begin{algorithm}[H]
    \begin{center}
        \begin{minipage}{1\linewidth}
            \footnotesize
            \caption{\texttt{FMG$(1,0)$-Cycle (\fmg) with }$D=I$}
            \label{alg:fmg}
            \begin{algorithmic}[1]
                \Require $A$, $b$, $P$, $R$, $\ell \ge 1$ (number of {\fmg } levels),\newline
                 $N \ge 1$ (number of {\ir } cycles with one V$(1,0)$ each), $m_A$, $m_P$, $m_R$

                \State $x \leftarrow 0$
                \COMMENT{Initialize \fmg}

                \If {$\ell > 1$}  \COMMENT{Check for Coarser Level}
                    \State $x_{\ell - 1} \leftarrow $\fmg$(A_{\ell - 1}, b_{\ell - 1}, P_{\ell - 1}, R_{\ell - 1}, \ell - 1, N, m_A, m_P, m_R)$\label{alg:fmg-interpolation}\\
                    \COMMENT{Compute Coarse-Level Approximation}

                    \State $P \leftarrow \texttt{quant}(P)$
                    \COMMENT{Quantize to $\wewe_\ell$ bits}

                    \State $x \leftarrow \texttt{qspmv}(P, x_{\ell - 1}, m_P, \wewe_\ell)$
                    \COMMENT{Interpolate Approximation $x \leftarrow P x_{\ell - 1}$}\label{fmg-interp}
                \EndIf
                \State $i \leftarrow 0$\COMMENT{Initialize {\ir } }
                \While {$i < N$}
                    \State {$x \leftarrow $\ir$(A, b, x, -1, m_A, \mathcal{V})$}
                    \COMMENT{Compute Correction by \ir-\V}

                    \State $i \leftarrow i + 1$
                    \COMMENT{Increment {\ir-\V} Cycle Counter}
                \EndWhile
                \State\Return $x$
                \COMMENT{Return Approximate Solution of $Ax = b$}
            \end{algorithmic}
        \end{minipage}
    \end{center}
\end{algorithm}

\end{document}